\documentclass[11pt]{article}

\usepackage{amssymb}

\newtheorem{theorem}{Theorem}

\newtheorem{corollary}[theorem]{Corollary}
\newtheorem{definition}[theorem]{Definition}
\newtheorem{examples}{Examples}
\newtheorem{lemma}[theorem]{Lemma}
\newtheorem{notation}{Notation}
\newtheorem{prop}[theorem]{Proposition}
\newtheorem{remark}{Remark}
\newtheorem{remarks}{Remarks}

\def\elaw{\stackrel{d}{=}}
\def\ss{{\,\succeq\,}}

\def\aa{{\eta}}
\def\CC{{\cal C}}

\def\ee{\varepsilon}
\def\esp{{\mathbb{E}}}
\def\FF{{\cal F}}
\def\F{{\cal G}}
\def\H{{\cal H}}
\def\II{{\cal I}}
\def\inte{{\mathbb{N}}}
\def\k{{\cal K}}
\def\lacc{\left\{}
\def\lcr{\left[}
\def\Lde{{\bf \tilde L}}
\def\LL{{\bf L}}
\def\lpa{\left(}
\def\lva{\left|}
\def\MC{{\cal M}}
\def\Nde{{\bf \tilde N}}
\def\NN{{\bf N}}
\def\pb{{\mathbb{P}}}
\def\R{{\mathbb{R}}}
\def\racc{\right\}}
\def\rcr{\right]}
\def\Rde{{\tilde R}}
\def\RC{{\cal R}}
\def\rpa{\right)}
\def\rva{\right|}
\def\tk{t^0_{jn}}
\def\tkk{t^1_{jn}}
\def\tkkk{t^2_{jn}}
\def\tkl{t^{-}_{jn}}
\def\Un{{\bf 1}}
\def\UU{{\bf U}}
\def\W{M}
\def\z{{\mathbb Z}}
\def\Zde{{\tilde Z}}

\newcommand{\dsty}{\displaystyle}
\newcommand{\fin}{\vspace{-0.3cm}
                  \begin{flushright}
                  \mbox{$\Box$}
                  \end{flushright}
                  \noindent}

\begin{document}

\begin{center}

\huge
{\bf Small deviations for fractional stable processes}

\vspace{8mm}

\Large
{\bf Mikhail Lifshits and Thomas Simon}

\end{center}

\vspace{2mm}

\begin{abstract}

\vspace{2mm}

\noindent
{Let $\{R_t, \; 0\leq t\leq 1\}$ be a symmetric $\alpha$-stable
Riemann-Liouville process with Hurst parameter $H > 0$.
Consider a translation invariant, $\beta$-self-similar, and
$p$-pseudo-additive functional semi-norm  $||.||$. We show that if
$H>\beta +1/p$ and $\gamma = (H-\beta-1/p)^{-1}$, then
$$\lim_{\ee\downarrow 0}\ee^{\gamma}\log\pb\lcr ||R||\leq \ee \rcr
 \; =\; - \k\;\in\; [-\infty, 0),$$
with $\k$ finite in the Gaussian case $\alpha = 2$.
If $\alpha < 2$, we prove that $\k$ is finite when
$R$ is continuous and $H>\beta +1/p + 1/\alpha$.
We also show that under the above assumptions,
$$\lim_{\ee\downarrow 0}\ee^{\gamma}\log\pb\lcr ||X||\leq \ee \rcr
 \; =\; - \k\;\in\; (-\infty, 0),$$
where $X$ is the linear $\alpha$-stable fractional motion
with Hurst parameter $H\in(0,1)$
(if $\alpha = 2$, then $X$ is the classical fractional
Brownian motion).
These general results cover many cases previously studied in the
literature, and also prove the existence of new small deviation constants,
both in Gaussian and Non-Gaussian frameworks.}
\end{abstract}

\vspace{2mm}

\def\abstractname{R\'esum\'e}

\begin{abstract}

\vspace{2mm}

\noindent
{Soit $\{R_t, \; 0\leq t\leq 1\}$ un processus de Riemann-Liouville
$\alpha$-stable sym\'etrique avec param\`etre de Hurst $H > 0$.
Consid\'erons une semi-norme fonctionnelle $||.||$ invariante par
translation, $\beta$-autosimilaire et $p$-pseudo-additive.
Nous montrons que si $H>\beta +1/p$ et $\gamma = (H-\beta-1/p)^{-1}$
alors
$$\lim_{\ee\downarrow 0}\ee^{\gamma}\log\pb\lcr ||R||\leq \ee \rcr
 \; =\; - \k\;\in\; [-\infty, 0),$$
avec $\k$ finie dans le cas gaussien $\alpha = 2$.
Lorsque $\alpha < 2$, nous montrons que $\k$ est finie quand $R$ est
continu et $H>\beta +1/p + 1/\alpha$.
Nous montrons aussi que sous ces hypoth\`eses
$$\lim_{\ee\downarrow 0}\ee^{\gamma}\log\pb\lcr ||X||\leq \ee \rcr
 \; =\; - \k\;\in\; (-\infty, 0),$$
o\`u $X$ est le mouvement fractionnaire lin\'eaire $\alpha$-stable
avec param\`etre de Hurst $H\in(0,1)$ (lorsque $\alpha = 2$, $X$
est le mouvement brownien fractionnaire usuel). Ces r\'esultats
g\'en\'eraux recouvrent de nombreux cas pr\'ec\'edemment \'etudi\'es
dans la litt\'erature et prouvent l'existence de nouvelles constantes
de petites d\'eviations, aussi bien dans le cadre gaussien que
non gaussien.}
\end{abstract}

\vspace{4mm}


\vspace{2mm}


\section{Introduction}

Let $X$ be a random process whose sample paths belong to some
functional normed
space $\lpa \FF, \|\cdot\|\rpa$. The investigation of the
small deviations (or small ball probabilities) of $X$ deals with the
asymptotics of
$$\pb\lcr ||X||\leq \ee \rcr$$
when $\ee\downarrow 0$, and has proved to be a difficult problem with
increasing
number of applications in Probability, Analysis, Complexity ... etc.
We refer to
the recent surveys \cite{LS} \cite{Lf} for a detailed account on this
subject.

In the literature, this problem is usually studied for a particular
class of processes and under a particular norm. It remains a
great challenge, a kind of "mission impossible", to find some 
principle describing small deviations for general classes
of processes and norms, rather than tackle the problem case by case.

The unique  successful attempt in this
direction was made by W. Stolz \cite{St1} \cite{St2},
who obtained estimates
\begin{equation}
\label{estim}
 -\infty \; < \; \liminf_{\ee\downarrow 0}\ee^{\gamma}\log\pb\lcr
||X||\leq \ee \rcr \; \leq \; \limsup_{\ee\downarrow 0}
\ee^{\gamma}\log\pb\lcr ||X||\leq \ee \rcr \; < \; 0,
\end{equation}
where $\{X_t, \; 0\leq t\leq 1\}$ is a Brownian motion (or more generally
a continuous Gaussian process with covariance function similar to that of
fractional Brownian motion), and $\gamma$ a finite positive parameter
 depending on the behavior of the (semi-)norm $\|\cdot\|$ on
linear combinations of Schauder functions. Many classical semi-norms fell
 into Stolz's scope: $L_p$-norms, H\"older and Sobolev semi-norms,
Besov norms... etc, and his estimates resumed from a general point of
view a lot of situations which were previously studied (see the references
 quoted in \cite{St1} \cite{St2}).

The next important issue is the {\em existence of the limit}
in (\ref{estim}), and this is the matter of the present paper.
Our main result says that if $||.||$ is a translation invariant,
$\beta$-self-similar and $p$-superadditive functional semi-norm -
see Definition 1 for more precisions about these notions, and
$\{R_t, \; 0\leq t\leq 1\}$ is a symmetric $\alpha$-stable
Riemann-Liouville process with Hurst parameter $H > \beta + 1/p$
($R$ can be viewed as a fractionally integrated symmetric
$\alpha$-stable L\'evy process, see Section 2 for a precise definition),
 then
\begin{equation}
\label{limit1}
\lim_{\ee\downarrow 0}\ee^{\gamma}\log\pb\lcr ||R||\leq \ee \rcr
 \; =\; - \k\;\in\; [-\infty, 0)
\end{equation}
with the rate $\gamma = (H-\beta-1/p)^{-1}$. From the technical point of view, 
the main ingredient of the
proof is a stochastic superadditive inequality which
is based upon the extrapolation-homogeneity of $R$, and
is then easy to combine with an exponential Tauberian theorem
and the standard subadditivity arguments.

Our framework has two secondary, but non-negligible advantages with respect to \cite{St2}. {\em Non-Gaussian} stable processes are included, as well as smooth Gaussian ones like integrated Brownian motion of arbitrary parameter.

A special effort is needed to prove that the constant $\k$
in (\ref{limit1}) is {\em finite}. In order to get the suitable
lower bound on the small deviation probabilities, we use roughly
the same method as in \cite{St1} \cite{St2}, except that we decompose
$R$ along Daubechies' wavelet bases, since Schauder system not
smooth enough when $H$ becomes too large.

In the last part of this paper, we extend the above results to
 a class of self-similar processes with long range dependence,
the so-called unilateral linear fractional stable motions.
These processes can be viewed as a possible generalization of
fractional Brownian motion with an underlying stable noise.
Thanks to a nice argument essentially due to
of Li and Linde \cite{LL1}, the problem is reduced
to a study of the Schauder decomposition of the associated
long memory process.

We conclude this article with a brief survey of concrete results.
It seems that our theorems synthesize everything that is known about the
existence of finite small deviation constants for continuous fractional processes
under translation-invariant semi-norms. Last but not least, we can harvest a bunch of new constants, both in Gaussian and Non-Gaussian situations.

Still, a major drawback of the above wavelet methods is that they
exclude discontinuous processes. As a rule, proving lower bound
probabilities for processes with jumps requires completely different
discretization techniques \cite{Mo} \cite{Ta}, and we have
no idea how much time it will take to enclose them into a global
strategy.

\section{Preliminaries} \label{norms}

\subsection{Parametrization: $(\beta, p)$-semi-norms}
Let $\II$ be the set of all closed bounded intervals of $\R$. Consider $\FF$ a linear space of functions from $\R$ to $\R$ and, for each $I\in\II$, let
$\FF_I$ be a linear space of functions from $I$ to $\R$ such that $f_I\in\FF_I$ for every $f\in\FF$, where $f_I$ stands for the restriction of $f$ to $I$.

We define a semi-norm $||\cdot||$ on $\FF$
as a family $\lacc ||\cdot||^{}_I, \; I\in\II\racc$
of functionals mapping $\FF_I$ to $\R^+$ such that
$||\lambda f||^{}_I = |\lambda| ||f||^{}_I$ and
$||f + g||^{}_I\leq ||f||^{}_I + ||g||^{}_I$
for every $\lambda\in\R, f, g\in\FF_I$.
We will use the notation $||f||^{}_I = ||f_I||^{}_I$ for every $f\in\FF$, $I\in\II$.
In the remainder of this paper we will assume that
$||\cdot||$ satisfies the following assumptions, which are verified by all
the classical semi-norms:

\vspace{2mm}

(A) $||\cdot||^{}_I\leq ||\cdot||^{}_J$ for every
$I, J\in\II$ such that $I\subset J$. (Contractivity)

\vspace{2mm}

(B) $||f||^{}_{I-c} = ||f(\cdot-c)||^{}_I$
for every $f\in\FF$, $I\in\II$ and $c\in\R$.
(Translation-invariance)

\begin{definition} Let $\beta\in\R$, $p\in (0, +\infty]$
and $||\cdot||$ be a contractive and translation-invariant
semi-norm on $\FF$. We say that $||\cdot||$ is an
{\em upper $(\beta, p)$-semi-norm} if it satisfies the
following properties:

\vspace{2mm}

{\em (C)} $||f(c\,\cdot)||^{}_{I/c} = c^\beta \ ||f||^{}_I$
for every $f\in\FF$, $I\in\II$ and
$c > 0.\;\;\;\;\;\;\;\;\;$ {\em ($\beta$-self-similarity)}

\vspace{2mm}

{\em (D)} For every $a_0<\ldots<a_n\in\R$ and $f\in\FF$

$$\lacc
\begin{array}{llr}
||f||^{}_{[a_0,a_n]} \;\geq\; {\lpa ||f||_{[a_0, a_1]}^p
+ \cdots + ||f||_{[a_{n-1},a_n]}^p\rpa}^{1/p}
& \mbox{{\em if} $p < +\infty$,} & \\

& & \;\mbox{\em ($p$-superadditivity)} \\

||f||^{}_{[a_0,a_n]} \;\geq\; \sup\lpa ||f||^{}_{[a_0, a_1]},
\ldots, ||f||^{}_{[a_{n-1},a_n]}\rpa
& \mbox{{\em if} $p = +\infty$.} &
\end{array}
\right.$$
\end{definition}

In the following we will denote by
$\UU(\beta,p)$ the set of upper $(\beta, p)$-semi-norms,
and set $\UU$ for the union of all $\UU(\beta,p)$'s.

\begin{remarks}{\em (a) Of course it suffices to take $n=2$
in the definition of $p$-superadditivity. We wrote the
property in this form in order to make it symmetric with the
$p$-subadditivity and the corresponding lower $(\beta, p)$-semi-norms,
which will appear just below.

\vspace{2mm}

\noindent
(b) From the inequality ${(a+b)}^q \geq a^q + b^q$
for every $a, b \geq 0$ and $q\geq 1$,
it follows that $\UU(\beta,p)\subset\UU(\beta,p')$
whenever $p'\geq p$.

\vspace{2mm}

\noindent
(c) In the definition of upper $(\beta,p)$-semi-norm, one can ask
for the possible values of the parameters $\beta$
and $p$.
By contractivity, homogeneity and translation-invariance, it is easy to see
that if
$ 0 < ||\Un||_{[0,1]} < \infty$ (where $\Un$ stands for
the unit function), then necessarily $\beta + 1/p \leq 0$. However,
the examples (c)-(f) below show that this inequality is no
more true whenever $||\Un||_{[0,1]} = 0$.}

\end{remarks}

We stress that most of the usual semi-norms belong to $\UU$:

\begin{examples}
{\em (a) The supremum semi-norm, which is given by
$$||f||^{}_I\; =\; \sup_{t\in I}|f(t)|$$
for every $I\in\II$, belongs to $\UU(0,+\infty)$.
\vspace{2mm}

\noindent
(b) The $L_p$-semi-norm, $p\geq 1$, which is given by
$$||f||^{}_I\; =\; {\lpa\int_I {\lva f(s)\rva}^p\, ds \rpa}^{1/p}$$
for every $I\in\II$, belongs to $\UU(-1/p,p)$.
\vspace{2mm}

\noindent
(c) The $\aa$-H\"older semi-norm, $0\leq \aa < 1$,
which is given by
$$||f||^{}_I\; =\; \sup_{s<t\in I}\frac{|f(t)-f(s)|}{{|t-s|}^{\aa}}$$
for every $I\in\II$, belongs to
$\UU(\aa,+\infty)$. In particular, the oscillation
semi-norm ($\aa = 0$) belongs to $\UU(0,+\infty)$. Similarly,
the Calder\'on-Zygmund semi-norm
$$||f||^{}_I\; =\; \sup_{s<t\in I}\frac{|2f((t+s)/2)-f(s) - f(t)|}{(t-s)}$$
belongs to $\UU(1,+\infty)$, and the $\aa$-Lipschitz semi-norm,
$\aa > 1$, which is given by
$$||f||^{}_I = \;
\sup_{s<t\in I}\frac{|f^{(n)}(t)-f^{(n)}(s)|}{{|t-s|}^{\aa-n}}$$
where $n < \aa < n+1$ (and by its Calder\'on-Zygmund
analogue for $\aa = n$), belongs to $\UU(\aa,+\infty)$.

\vspace{2mm}

\noindent
(d) The (strong) $p$-variation semi-norm, $p\geq 1$,
which is given by
$$||f||^{}_I\; =\; {\lpa \sup_{t_0 < \ldots < t_n \in I }
\sum_{i = 1}^n {|f(t_i)-f(t_{i-1})|}^p\rpa}^{1/p}$$
for every $I\in\II$, belongs to $\UU(0,p)$.

\vspace{2mm}

\noindent
(e) The $(\aa, p)$-Sobolev semi-norm, $p\geq 1$ and
$0\leq\aa +1/p < 1$, which is given by
$$||f||^{}_I\; =\; {\lpa\int_{I}\int_{I}
{\lpa\frac{|f(t)-f(s)|}
{{|t-s|}^{\aa +1/p}}\rpa}^p ds\, dt\rpa}^{1/p}$$
for every $I\in\II$, belongs to
$\UU(\aa -1/p, p)$.

\vspace{2mm}

\noindent
(f) The $(\aa, p, q)$-Besov semi-norm,
$\aa > 0$ and $p, q \geq 1$, which is given by
$$||f||^{}_I\; =\;
{\lpa\int_0^{|I|}{\lpa\frac{\omega_{p,I}(t,f)}
{t^{\aa}}\rpa}^q \frac{dt}{t}\rpa}^{1/q}$$
for every $I\in\II$, where
$$\omega_{p,I}(t,f) = \sup_{|h|\leq t}{\lpa\int_{I_h}
{|f(x-h)-f(x)|}^p dx\rpa}^{1/p},$$
$I_h = \lacc x\in I\,|\, x-h\in I\racc$ and $|I|$ stands
for the Lebesgue measure of $I$, belongs to $\UU(\aa-1/p, +\infty)$.
Notice that the usual $(\aa, p, q)$-Besov semi-norm,
 where the $L_p$-semi-norm is added in the definition
 (see e.g. \cite{CKR}), does {\em not} belong to $\UU$
since $L_p$-semi-norms have a different self-similarity index.
}

\end{examples}

It follows from Remark (b) above that the concept of upper
$(\beta, p)$-semi-norm is not sharp enough for our further purposes.
We need some kind of inverse property to $(D)$, in order to
work with a family of {\em disjoint} sets of semi-norms.
Let $\CC^l_K$ (resp. $\CC^0_K$) denote
the class of all $l$-times continuously differentiable (resp.
continuous) functions with compact support.
For technical reasons, we may have to make the following assumptions on $||\cdot||$ which, again, are verified by all the classical semi-norms:

\vspace{2mm}

(E) $\exists\;l\in\inte$ such that $\CC^l_K\subset \FF$.
\qquad (Smooth-finiteness)

\vspace{2mm}

(F) For every $I\in\II$ and $f, f_n \in\FF_I$,
$$f_n\rightarrow f\; \mbox{uniformly on $I$}\;\Longrightarrow\;||f||_I\leq\limsup ||f_n||_I.\;\;\;\;\;\;\;\mbox{(Lower semi-continuity)}$$

\begin{remark}{\em
Clearly, $l$-smooth-finiteness yields $||f||_I < +\infty$
for every $f$ of class $\CC^l$ and every $I\in\II$.}
\end{remark}

\begin{definition}
Let $\beta\in\R$, $p\in (0, +\infty]$
and $||\cdot||$ be a contractive, translation-invariant, smooth-finite,
lower semi-continuous and $\beta$-self-similar semi-norm on $\FF$.
We say that $||\cdot||$ is a {\em lower $(\beta, p)$-semi-norm}
if it satisfies the following property:

\vspace{2mm}

{\em (G)} There exists $C_p > 0$ such that for
 every $a_0<\ldots<a_n\in\R$ and $f\in\CC^l_K$ verifying
$f(a_0) = \cdots = f(a_n) = 0$,

$$\lacc\begin{array}{llr}
||f||^{}_{[a_0,a_n]} \;\leq\;
C_p {\lpa ||f||_{[a_0, a_1]}^p + \cdots +
||f||_{[a_{n-1},a_n]}^p\rpa}^{1/p} & \mbox{{\em if}
$p < +\infty$,} & \\

& & \mbox{\em ($p$-subadditivity)} \\

||f||^{}_{[a_0,a_n]} \;\leq\; C_{\infty}
\sup\lpa ||f||^{}_{[a_0, a_1]}, \ldots,
||f||^{}_{[a_{n-1},a_n]}\rpa & \mbox{{\em if} $p = +\infty$.}
 & \end{array}
\right.$$

\end{definition}

In the following we will denote by $\LL(\beta,p)$
the set of lower $(\beta, p)$-semi-norms, and set $\LL$
for the union of the $\LL(\beta,p)$'s.
We will say that $||\cdot||$
is a {\em $(\beta, p)$-semi-norm} if it belongs to
$\NN(\beta,p) = \UU(\beta,p)\cap\LL(\beta,p)$ and analogously
we will set $\NN$ for the union of $\NN(\beta,p)$'s.

\begin{remarks}{\em (a) Since the constant $C_p$ may be larger
than 1, it is not enough to take $n = 2$ in the
definition of $p$-subadditivity, contrary to $p$-superadditivity.

\vspace{2mm}

\noindent
(b) From the inequality
${(a_0 + \cdots + a_n)}^q \leq a_0^q + \cdots + a_n^q$
 for every $a_0,\ldots, a_n\geq 0$ and $0 \leq q\leq 1$,
 it follows that $\LL(\beta,p)\subset\LL(\beta,p')$ whenever $p'\leq p$.

\newpage
\vspace{2mm}

\noindent
(c) Again, one can ask for the possible values of the parameters $\beta$
 and $p$ in the definition of a lower $(p,\beta)$-semi-norm.
By contractivity, homogeneity and translation-invariance,
it is easy to see that when the inequalities in (G) hold for
the function $f=\Un$ and for all partitions, then necessarily $\beta + 1/p
\geq 0$. This inequality is actually true for all available examples,
 but we got stuck in proving that in full generality.

\vspace{2mm}

\noindent
(d) It is not difficult to see that $\NN(\beta, p)\cap\NN(\beta', p') =
 \emptyset$ as soon as $(\beta, p)\neq(\beta', p')$.}
\end{remarks}

Notice that each of the above examples (a)-(d) belongs to
$\LL$, and hence to $\NN$. Actually in each case we had chosen
the smallest possible parameter $p$, i.e. we could have written
$||\cdot||\in\LL(\beta, p)$ as well as $||\cdot||\in\UU(\beta, p)$.
This fact is trivial for examples (a)-(c) where we can take $C_p = 1$,
 a bit more involved in example (d) where we have to take
$C_p = 2^{1- 1/p}$ as soon as $p \geq 1$
(notice that here the condition $f(a_0) = \cdots = f(a_n) = 0$
 is essential).

However, $(\aa, p, q)$-Besov and $(\aa, p)$-Sobolev
semi-norms do not belong to $\LL$, because in these cases
the constant $C_p$ may depend on the subdivision  $a_0 < \ldots < a_n$.
 These two examples are important in certain contexts,
 and this is why we would like to introduce a weaker definition, which
will be given in terms of the evaluation of $||\cdot||$ along
specific families of functions.

If $l\in\inte$ and $\psi\in\CC^l_K$, we introduce
$\Psi=\lacc\psi_{jn}, \; n \in \z, \; j\ge 0 \racc$,
the two-parametric subset of $\CC^l_K$ defined by
$$\psi_{jn}(t)= \psi(2^j(t-n))$$
for every $n\in\z, j\ge 0, t\in\R$.

\begin{definition} Let $\beta\in\R$, $p\in(0,+\infty]$ and $||\cdot||$
 be a lower semi-continuous, $l$-smooth-finite semi-norm on $\FF$.
We say that $||\cdot||$ is a {\em lower $(\beta, p)$-semi-norm
in the wide sense} with respect to $\Psi$,
if it satisfies the following property:
\vspace{2mm}

{\em ($\rm \tilde G$)} There exists $C_p > 0$ such that for every
$j\ge 0$, arbitrary $x_1,\ldots,x_m\in\R$, and arbitrary
$n_1,\ldots, n_m\in\inte$
such that the supports of the functions $\{\psi_{j,n_i}, \; 1\le i\le m\}$
have disjoint interiors,
$$\lacc
\begin{array}{ll}
{\lva\lva {\dsty
 \sum_{i=1}^{m} x_i \psi_{jn_i} }\rva\rva}_{[0,1]} \;
\leq\; C_p\, 2^{\beta j}   \
{\lpa {\dsty   \sum_{i=1}^{m} |x_i|^p }\rpa}^{1/p}
 & \mbox{{\em if} $p < +\infty$,}  \\

&  \\

{\lva\lva {\dsty \sum_{i=1}^{m} x_i \psi_{jn_i}  }\rva\rva}_{[0,1]}  \;
\leq\; C_\infty\, 2^{\beta j}  \
 {\dsty   \sup_{1\le i\le m} |x_i|  }
& \mbox{{\em if} $p = +\infty$.}
\end{array}
\right.$$

\end{definition}

Similarly we will denote by $\Lde(\beta,p)$ the set of lower
$(\beta, p)$-semi-norms in the wide sense, and set $\Lde$
for the union of the $\Lde(\beta,p)$'s. Analogously, we define
$\Nde(\beta,p) = \UU(\beta,p)\cap\Lde(\beta,p)$ and set $\Nde$
for the union of the $\Nde(\beta,p)$'s.

\vspace{2mm}

\begin{remarks}{\em (a) In the sequel, $\Psi$ will be either a family
of sufficiently smooth wavelet functions, or some Schauder system on $[0,1]$.

\vspace{2mm}

\noindent
(b) Condition ($\rm \tilde G$) just means
that (G) holds for linear combinations of specific functions
with the same "frequency", and it is a well-known condition to
obtain lower bounds for small deviation probabilities in a Gaussian
framework \cite{St1} \cite{St2}. Actually we will use it
for the same purposes, but sometimes in a more general context, see
Section \ref{lower}.

\vspace{2mm}

\noindent
(c) At first sight, neither translation-invariance
 nor self-similarity are required in the definition of $\Lde(\beta,p)$.
Actually, working with a specific family of function allows us
to combine these necessary properties in the single inequality given
by ($\rm \tilde G$).

\vspace{2mm}

\noindent
(d) Again we can prove that $\LL(\beta,p)\subset\LL(\beta,p')$
(resp. $\Nde(\beta, p)\cap\Nde(\beta', p') = \emptyset$)
whenever $p'\leq p$ (resp. $(\beta, p)\neq(\beta', p')$).

\vspace{2mm}

\noindent
(e) It is a bit tedious but not difficult to see that
$(\aa, p, q)$-Besov semi-norm belongs to $\Nde(\aa,+\infty)$
 and that  $(\aa, p)$-Sobolev semi-norm belongs to
$\Nde(\aa -1/p,p)$, both with respect to the Schauder system.
}
\end{remarks}

\begin{notation} {\em In the following we will mainly
consider functions restricted to $[0,1]$. By minor abuse of notations,
and for the sake of brevity, we will set $||f|| = ||f||_{[0,1]}$
for every $f\in\FF$, although such statements like
$"||\cdot||\in\UU(\beta, p)"$ or $"||\cdot||\in\LL(\beta, p)"$
 will always refer to the family
$||\cdot|| = \lacc ||\cdot||^{}_I, \; I\in\II\racc$.
By $C$ we will always mean a positive finite constant
independent of the involved parameters, and
whose value may change from line to line.}
\end{notation}

\subsection{Riemann-Liouville processes and their
associated linear fractional stable motions}

Let $\lacc Z_t, \; t\geq 0\racc$ be a symmetric
$\alpha$-stable
process with index $\alpha\in(0,2]$, i.e. $Z$ is a
L\'evy process  whose characteristic function is given by
$$\esp\lcr e^{\rm{i} \lambda Z_t}\rcr\; =\; e^{-t|\lambda|^{\alpha}}$$
for every $t\geq 0$ and $\lambda\in\R$.
It is well-known and easy to see that for every $H > 0$ the
following stochastic integral:
$$R^H_t\; =\;\int_0^t (t-s)^{H-1/\alpha} dZ_s,$$
is well-defined for every $t > 0$. We set $R^H_0 = 0$ and
call $\lacc R^H_t, \; t\geq 0\racc$ the {\em Riemann-Liouville process}
 (in abridged form: RLP)
with Hurst parameter $H$. This latter terminology is motivated
 by the following $H$-self-similarity property of $R^H$: for every $c > 0$
$$\lacc R^H_{ct}, \; t\geq 0\racc\;\elaw\;\lacc c^H R^H_t, \;
t\geq 0\racc.$$
When no confusion is possible, we will drop the subscript
$H$ and write $R=R^H$ for the sake of brevity. Notice that $R$
has no stationary increments, unless $H = 1/\alpha$.
However $R$ shares some kind of extrapolation-homogeneity,
which will be important in the sequel: namely, if we set
$$R_{a,t}\; =\;\int_a^t (t-s)^{H-1/\alpha} dZ_s$$
for every $t\geq a\geq 0$, then the following equality in
law holds:
$$\lacc R_{a,a+t}, \; t\geq 0\racc\;\elaw\;\lacc R_t, \; t\geq 0\racc.$$
Besides, $R$ shares an equally important independence property:
 if $\Rde_{a,t} = R_t - R_{a,t}$, then for every $a \geq 0$ the processes
 $\lacc R_{a,a+t}, \; t\geq 0\racc$ and $\{\Rde_{a,a+t}, \; t\geq 0\}$
are independent. These three properties are easy to see
as respective direct consequences (in reverse order) of the independence,
 stationarity, and stability of the increments of $Z$.

Of course, $Z$ itself is an RLP with Hurst parameter $H = 1/\alpha$.
 Notice also that up to normalization constants,
the family of Riemann-Liouville processes is closed with respect to
time-integration. In particular, the $m$-times integrated Brownian motion
is an RLP with parameters $\alpha=2$ and $H=m+ 1/2$.

\vspace{2mm}

The Riemann-Liouville process is closely related to
$\lacc X^H_t,\; t\geq 0\racc$, the so-called
{\em linear stable fractional motion} (in abridged form: LFSM)
with Hurst parameter $H$. $X^H$ can be defined through the following
(independent) decomposition: $X^H = R^H + \W^H$ where $\W^H_0 = 0$ and
$$\W^H_t=\int_0^{+\infty} \lpa (t+s)^{H-1/\alpha} - s^{H-1/\alpha}
\rpa\; d\Zde_s$$
for every $t> 0$, $\Zde$ being an independent copy of $Z$.
We call $\W^H$ the {\em long memory process} (in abridged form: LMP)
associated to $X^H$. Notice that the stochastic integral defining
$\W^H$ diverges a.s. at $+\infty$ as soon as $H\geq 1$. Hence,
$X^H$ can be defined only for $H\in(0,1)$.
Again, when no confusion is possible, we will drop the
subscript $H$ and write $X = X^H$ (resp. $\W = \W^H$) for
the sake of brevity.

$X$ is also $H$-self-similar but, contrary to $R$, does
not share the above extrapolation-homogeneity property.
Instead, its increments are truly stationary:
$$\lacc X_{s+t} - X_s, \; t\geq 0\racc\;\elaw\;\lacc X_t, \;
t\geq 0\racc$$
for every $s \geq 0$, and $X$ is a so-called {\em $H$-sssi process}.
We refer to the monographs \cite{ST} \cite{EM} for an extensive account on
these latter processes. Notice that in the Gaussian case $\alpha = 2$,
$X$ coincides with the canonical fractional Brownian motion, up to
some numerical factor (see Proposition 7.2.6 in \cite{ST}).
Moreover, an alternative definition of $X$ can
then be given through
the following ("well-balanced") decomposition:
$$X_t=\int_{-\infty}^{+\infty} \lpa |t+s|^{H-1/2} -
|s|^{H-1/2} \rpa\; \dot{B}(ds),$$
where $\dot{B}$ is the usual white noise
(see Exercise 7.2 in \cite{ST}). When $\alpha < 2$,
the corresponding integral
\begin{equation} \label{balanced}
\tilde X_t=\int_{-\infty}^{+\infty} \lpa |t+s|^{H-1/\alpha} -
|s|^{H-1/\alpha} \rpa\;\dot{Z}(ds)
\end{equation}
also makes sense, where $\dot{Z}$ is the so-called symmetric $\alpha$-stable noise. However, $\tilde X$ is no more equivalent to $X$
(see Theorem 7.4.5 in \cite{ST}). Actually, many other "bilateral"
 $H$-sssi $\alpha$-stable processes can be constructed,
which are all non-equivalent except in the Gaussian case
(see Definition 7.4.1 and Theorem 7.4.5 in \cite{ST}).
 Hence, our unilateral LFSM process is just {\em one} possible
 stable generalization of fractional Brownian motion.
 However, we will restrict our study to this unilateral
 process, even though our results can probably be adapted to
 some other situations.

Notice that the process $M$ is {\em smooth} on $(0, +\infty)$,
so that $X$ and $R$ exhibit similar local properties. In particular it
is well-known (see Chapters 9-12 in \cite{ST}) that $R$ (resp. $X$)
admits a continuous version if and only if
$$\alpha = 2\;\;\;\;\;\mbox{or}\;\;\;\;\; H > 1/\alpha.$$
In Section 6 we will give some other local properties of $R$ and $X$,
related to the examples of semi-norms (a)-(f) listed above. We finally
refer to \cite{MaR} for a more thorough comparison between $R$ and $X$,
in the context of econometric applications.

\section{Existence of the small deviation constant for RLP}
\label{main}

We can now state the main result of this paper.

\begin{theorem} \label{t31} Let $||\cdot||\in\UU(\beta, p)$ and $R$
be an RLP with Hurst parameter $H$. Assume that $H>\beta +1/p$
and  set $\gamma = (H-\beta-1/p)^{-1}$.
Then there exists $\k\;\in\; (0, +\infty]$ such that
\begin{equation}
\lim_{\ee\downarrow 0}\ee^{\gamma}\log\pb\lcr ||R||\leq \ee \rcr
 \; = \; - \k.
\end{equation}
\end{theorem}

\begin{remarks}{\em (a) It is interesting to note that the
 stability index $\alpha$ does not directly show up in the
 expression of the small deviation rate.

\vspace{2mm}

\noindent
(b) This result says nothing about the finiteness of the constant
 $\k$, which is of course a very important feature. To this aim,
one needs to show a lower bound for small deviation probabilities
with appropriate order - see Section \ref{lower}.
While {\it a priori}
nothing indicates that the order $\gamma$ of small deviation
probabilities we propose in Theorem \ref{t31} is the right one,
it is indeed correct for {\it all} available examples, as soon
as $p$ is chosen as small as possible - see Section \ref{survey}.
}
\end{remarks}

The next two paragraphs are devoted to the proof of
Theorem \ref{t31} and we fix once and for all
an upper $(\beta, p)$-semi-norm $||\cdot||$, and an RLP $R$
with Hurst parameter $H>\beta +1/p$. We begin in proving a kind of
stochastic superadditivity property, which will be our crucial
argument.

\subsection{Stochastic superadditivity}

Let $X$ and $Y$ be two real-valued random variables.
 We say that $X$ is {\em stochastically larger} than $Y$ if and only if
$$\pb\lacc X \le r\racc\; \leq \; \pb \lacc Y \leq r\racc$$
for every $r\in\R$, and we write in this case $X\; \ss\; Y$.
 If $X$ and $Y$ are positive, then clearly $X\; \ss\; Y$ entails that
$$\esp\lcr\exp -\lambda X\rcr\; \leq \;\esp\lcr\exp -\lambda Y\rcr$$
for every $\lambda \geq 0$.

\begin{prop} \label{in1} Suppose that $p < \infty$. Let $R_1$,
$R_2$ be two independent copies of $R$ and set $q = p(H-\beta) > 0$.
Then for every $a, b\geq 0$
$$ (a+b)^q \ ||R||^p\; \ss\; a^q\ ||R_1||^p + b^q \ ||R_2||^p.$$
\end{prop}

\noindent
{\em Proof.} Fix $a, b > 0$ and set $c = a + b$.
We first appeal to the $p$-superadditivity of $||\cdot||$ and get
$$||R||_{[0,c]}^p\; \ge\; ||R||_{[0,a]}^p + ||R||_{[a,c]}^p.
$$
By self-similarity for both $||\cdot||$ and $R$, we also have
\begin{eqnarray*}
 || R ||_{[0,c]}^p & \elaw & ||c^H R(c^{-1}\cdot)||_{[0,c]}^p \;
=\; c^{p(H -\beta)}|| R ||_{[0,1]}^p \; =\; c^q || R ||^p
\end{eqnarray*}
and, similarly,
$$ ||R||_{[0,a]}^p \;\elaw\;  a^q || R_1 ||^p\;\;\; \mbox{and}\;\;\;
||R||_{[0,b]}^p \;\elaw\;  b^q || R_2 ||^p.
$$
Putting everything together, we now see that it is enough to show
\begin{equation}
  \label{bp1}
||R||_{[0,a]}^p + ||R||_{[a,c]}^p\; \ss\;
||R||_{[0,a]}^p + ||R_a||_{[a,c]}^p,
\end{equation}
because of the translation-invariance of $||\cdot||$ and the
 extrapolation-homogeneity of $R$.
Here and throughout this section, $R_a$ denotes the process
$\{R_{a,t},\, t\ge a\}$
For every $a\geq 0$,
set $\F_a$ for the filtration generated by
$\lacc Z_s,\; 0 \leq s \leq a\racc$ and $\pb_a$
for the conditional probability with respect to $\F_a$.
We clearly have
\begin{eqnarray*}
\pb\lcr\ ||R||_{[0,a]}^p + ||R||_{[a,c]}^p \leq r \rcr & = &
 \esp\lcr\pb_a \lcr\ ||R||_{[a,c]}^p \leq r -  ||R||_{[0,a]}^p \rcr\rcr \\
& = & \esp \lcr \pb_a \lcr ||\tilde R_a + R_a||_{[a,c]}^p \leq r -
||R||_{[0,a]}^p \rcr\rcr
\end{eqnarray*}
for every $r\in\R$. Since $R_a$ is a conditionally Gaussian
process under $\pb_a$ (see \cite{ST} p.153-154), we can apply
Anderson's inequality under $\pb_a$ and get (recall that $\tilde R_a$ is
 $\F_a$-measurable)
$$ \pb_a\lcr ||\tilde R_a + R_a ||_{[a,c]}^p \le r -
||R||_{[0,a]}^p \rcr\; \leq\;\pb_a\lcr||R_a||_{[a,c]}^p \leq r -
||R||_{[0,a]}^p \rcr.$$
Now, since $R_a$ is independent of $\F_a$, we can average backwards
 and obtain
$$\pb\lcr ||R||_{[0,a]}^p + ||R||_{[a,c]}^p \leq r\rcr\; \leq\;\pb\lcr
||R||_{[0,a]}^p + ||R_a||_{[a,c]}^p \leq r \rcr$$
for every $r\in\R$, as desired for (\ref{bp1}).

\fin

The next proposition covers the simpler situation when $p=\infty$,
and follows mainly the outline of Theorem 2.1 in \cite{LL1}.

\begin{prop} \label{in2}
Suppose that $ p = +\infty$ and set $q= H - \beta > 0$. For every $a,
b\geq 0$ and $r\in\R$
$$
\pb\lcr (a+b)^q ||R||\,\leq\, r \rcr\; \leq\; \pb\lcr a^q ||R||\,\leq\,r
\rcr\;\pb\lcr b^q ||R||\,\leq\, r \rcr.
$$
\end{prop}

\noindent {\em Proof.} Fix $a, b > 0$ and set $c = a + b$. Again we appeal
 to the $\infty$-superadditivity of $||\cdot||$ and get
$$||R||_{[0,c]} \geq \sup\lacc ||R||_{[0,a]}, ||R||_{[a,c]} \racc.$$
Therefore,
\begin{eqnarray*}
\pb\lcr ||R||_{[0,c]}\leq r \racc & \leq & \pb\lcr ||R||_{[0,a]}\leq r,
 ||R||_{[a,c]}\leq r\racc \; =\; \esp\lcr {\bf 1}_{\lacc ||R||_{[0,a]
}\leq r \racc}
\pb_a \lcr ||R||_{[a,c]}\le r\rcr\rcr
\end{eqnarray*}
for every $r\in\R$, where $\pb_a$ is defined as above.
By Anderson's inequality, independence of $R_a$ with respect
to $\F_a$, extrapolation-homogeneity of $R$ and
translation-invariance of $||\cdot||$, we obtain
\begin{eqnarray*}
\pb_a \lcr ||R||_{[a,c]}\le r\racc & = & \pb_a \lcr
 ||\tilde R_a+R_a||_{[a,c]}\le r\rcr\;\leq\; \pb_a \lcr
||R_a||_{[a,c]}\le r\rcr \; =\;\pb \lcr||R||_{[0,b]}\le r\rcr.
\end{eqnarray*}
Averaging backwards, we get
$$\pb\lcr ||R||_{[0,c]}\le r \rcr\; \leq\;
\pb\lcr ||R||_{[0,a]}\le r\rcr\pb\lcr||R||_{[0,b]}\le r\rcr
$$
for every $r\in\R$. Now this yields the desired inequality
$$
\pb\lcr c^q ||R||\,\leq\, r \rcr\; \leq\; \pb\lcr a^q ||R||\,
\leq\,r \rcr\;\pb\lcr b^q ||R||\,\leq\, r \rcr,
$$
since by self-similarity of $R$ and $||\cdot||$
$$\pb\lcr ||R||_{[0,t]} \leq r \rcr\; =\; \pb\lcr t^q ||R|| \leq r \rcr$$
for every $t \geq 0$ and $r\in\R$.

\fin

We are now ready to proceed to the proof of Theorem \ref{t31}.

\subsection{Proof of Theorem \ref{t31}}

Again we first consider the case $p<\infty$. Proposition \ref{in1}
yields the following decisive inequality for Laplace transforms:
$$ \esp\lcr \exp - (a+b)^q ||R||^p \rcr \;\le\; \esp\lcr \exp -a^q
||R||^p \rcr \ \esp\lcr \exp - b^q ||R||^p \rcr
$$
for every $a, b \geq 0$ and $q = p(H-\beta)$ as above.
This entails that the log-Laplace transformation of $||R||$ defined by
$$\Psi(h) = \log \esp \lcr \exp -h^{q} ||R||^p \rcr$$
for every $h \geq 0$, is a continuous negative function which satisfies
$\Psi(a+b) \leq \Psi(a) + \Psi(b)$ for every $a, b \geq 0$.
 By the standard subadditivity argument, we obtain
$$\lim_{h\to\infty} \frac{\Psi(h)}h \; = \; \inf_{h \geq 0}
\frac{\Psi(h)}h\; = \; -C \;\in \; [-\infty, 0)$$
and, returning to the Laplace transform,
$$\lim_{\lambda\to +\infty}\lambda^{1/q}\; \esp\lcr \exp -\lambda
||R||^p \rcr\; =\; -C.$$
Notice that $q>1$ by assumption. Hence we can appeal to de Bruijn's
exponential Tauberian theorem (see Theorem 4.12.9 in \cite{BGT}, or
Theorem 3.5 in
\cite{LS} for a more comfortable formulation), which yields
$$\lim_{\ee\to 0}\ee^{1/(q-1)}\;\pb\lcr ||R||^p \le \ee \rcr\;
 =\; -\k \; =\; -(q-1){\lpa C/ q \rpa}^{q/(q-1)},$$
and finally
$$\lim_{\ee\to 0}\ee^{\gamma}\;\pb\lcr ||R||
\le \ee \rcr\; =\; - \k \;\in \; [-\infty, 0)$$
with $\gamma = (H-\beta-1/p)^{-1}$.
This completes the proof of Theorem \ref{t31} when $p < +\infty$.
\vspace{2mm}

In the case $p= +\infty$ we do not even need Laplace transformation,
since it follows directly from Proposition \ref{in2} that the function
$$ \Psi(h) = -\log \pb\lcr ||R|| \le h^{\beta-H} \rcr$$
is subadditive. Again this entails
$$\lim_{h\to\infty} \frac{\Psi(h)}h \; = \; \inf_{h \geq 0}
 \frac{\Psi(h)}h\; = \; -\k \;\in \; [-\infty, 0),$$
and we obtain
$$\lim_{\ee\to 0}\ee^{\gamma}\;\pb\lcr ||R|| \le \ee \rcr\; =\; -
 \k \;\in \; [-\infty, 0)$$
with $\gamma = (H-\beta)^{-1}$, as desired when $p = +\infty$.
\fin

\section{Lower bounds: finiteness of the constant for continuous
RLP} \label{lower}

In this section we obtain a suitable lower bound for small deviation
 probabilities which will allow us to prove, under certain conditions,
that the constant $\k$ from Theorem \ref{t31} is actually finite
 whenever $||\cdot||\in \LL(\beta, p)$ as well. Unfortunately,
our method is only efficient in the {\em continuous} case, i.e. when
$$\alpha = 2\;\;\;\;\;\mbox{or}\;\;\;\;\; H > 1/\alpha.$$
An explanation for this understandable, but important restriction on
$R$ will be given later. Our result reads as follows:

\begin{theorem} \label{t41} Let $||\cdot||\in\LL(\beta, p)$ and $R$ be a
 continuous $\alpha$-stable RLP with Hurst parameter $H$. Suppose that
$H > \beta + 1/p$ if $\alpha = 2$ and $H > \beta+ 1/p+1/\alpha$ if
 $\alpha < 2$. Then
\begin{equation} \label{linf}
\liminf_{\ee\downarrow 0}\ee^{\gamma}\log\pb\lcr
||R||\leq \ee \rcr
 \; > \; -\infty,
\end{equation}
with $\gamma = (H-\beta-1/p)^{-1}$.
\end{theorem}

Combining Theorem \ref{t31} and Theorem \ref{t41} yields readily the
following fairly general small deviation theorem for continuous
Riemann-Liouville processes:

\begin{corollary}  \label{c41} Let $||\cdot||\in\NN(\beta, p)$
and $R$ be a continuous $\alpha$-stable RLP with Hurst parameter $H$.
Suppose that $H > \beta + 1/p$ if $\alpha = 2$ and $H > \beta+ 1/p+1/\alpha$
 if $\alpha < 2$. Then there exists $\k\in (0,\infty)$ such that
$$\lim_{\ee\downarrow 0}\ee^{\gamma}\log\pb\lcr ||R||\leq \ee \rcr
 \; = \; -\k, $$
with $\gamma = (H-\beta-1/p)^{-1}$.
\end{corollary}

We stress that if the parameter $H$ is not too large,
more precisely if $H<2$, then it is possible to obtain the lower bound
for small deviation probabilities
under a weaker assumption than $||\cdot||\in\LL(\beta, p)$. Here we
just need that $||\cdot||\in\Lde(\beta, p)$ with respect to
$$\Psi\; =\; \lacc\psi_{jn}=\psi(2^jt-n+1),\; 1\le
n \le 2^j,\; j\ge 0\racc,$$
the Schauder system generated by the triangular
 function $\psi(t)=\Un_{[0,1]} (t)(1-|2t-1|)$.

Notice that on each level $j$ the supports of the functions $\psi_{jn}$
have disjoint interiors,
so that here $||\cdot||\in\Lde(\beta, p)$ w.r.t. $\Psi$ simply means that
$$ {\lva\lva
           {\dsty \sum_{n=1}^{2^j} x_n \psi_{jn} }
   \rva\rva} \;
\leq\; C_p\, 2^{\beta j} {\lpa {\dsty   \sum_{n=1}^{2^j}
|x_n|^p }\rpa}^{1/p}$$
for $p < +\infty$, and with an obvious modification for $p=\infty$.
We have an analogous result to Theorem \ref{t41}:
\begin{theorem} \label{t42} Let $||\cdot||\in\Lde(\beta, p)$ w.r.t.
Schauder system $R$ be a continuous $\alpha$-stable RLP with Hurst parameter
$H < 2$. Suppose that $H > \beta + 1/p$ if $\alpha = 2$ and $H >
\beta+ 1/p+1/\alpha$ if $\alpha < 2$. Then
\begin{eqnarray*}
\liminf_{\ee\downarrow 0}\ee^{\gamma}\log\pb\lcr
||R||\leq \ee \rcr
 & > & -\infty,
\end{eqnarray*}
with $\gamma = (H-\beta-1/p)^{-1}$.
\end{theorem}
And, of course, we get the corresponding corollary:
\begin{corollary} \label{c42}  Let $||\cdot||\in\Nde(\beta, p)$ w.r.t.
Schauder system and $R$ be a continuous $\alpha$-stable RLP
with Hurst parameter $H < 2$. Suppose that
$H > \beta + 1/p$ if $\alpha = 2$ and $H > \beta+ 1/p+1/\alpha$
 if $\alpha < 2$. Then there exists $\k\in (0,\infty)$ such that
$$\lim_{\ee\downarrow 0}\ee^{\gamma}\log\pb\lcr ||R||\leq \ee \rcr
 \; = \; -\k, $$
with $\gamma = (H-\beta-1/p)^{-1}$.
\end{corollary}

\subsection{Decorrelating stable arrays}

In this paragraph we prove a crucial lower bound, which is a
generalization of Lemma 2.1 in \cite{St2} to arrays of symmetric
$\alpha$-stable random variables, and which will be useful for
both Theorems \ref{t41} and \ref{t42}. It relies on a version of
\v Sid\'ak's inequality for stable variables recently obtained in
 this framework by G. Samorodnitsky \cite{Sa}.

\begin{lemma} \label{samo}
Let $M, h > 0$ and $\lacc y_{jn},\; 1\leq n \leq M 2^{h j}, \;
 j\geq 0\racc$ be an array of identically distributed
symmetric $\alpha$-stable
 random variables. Let $z, \delta > 0$ be such that
 $\delta< z$ if $\alpha = 2$ and $ \delta < z - h/\alpha$ if
$\alpha < 2$. Let $m>0$ be an integer and set
$$d_j\; =\; d_j(m)\; =\;
2^{z(j-m) - \delta|j-m|}. $$
Then there exists a constant $C$ depending only on $M,h, z,
 \delta$ such that
$$\pb\lcr |y_{jn}| \leq d_j, \; 1\leq n \leq M 2^{h j}, \; j\geq 0\rcr \;
 \geq \; \exp - C 2^{h m}.$$
\end{lemma}

\noindent{\em Proof}. Up to normalization, the case $\alpha = 2$ is
just the statement of Lemma 2.1 in \cite{St2}. Hence we can concentrate
 on the case $\alpha < 2$, and first notice that Lemma 2.1 in \cite{Sa}
 entails the following decorrelation inequality:
$$\pb\lcr |y_{jn}| \leq d_j, \; 1\leq n \leq M 2^{h j}, \; j\geq 0\rcr \;
\geq \; \prod_{j\geq 0}
{\pb\lcr |y| \leq d_j\rcr}^{M 2^{h j}},$$
where $y$ is some symmetric $\alpha$-stable
 random variable. We can decompose the right-hand side into
$$\prod_{j\geq 0} {\pb\lcr |y| \leq d_j\rcr}^{M 2^{h j}}\; =\;
 \lpa\prod_{j\geq m} {\pb\lcr |y| \leq d_j\rcr}^{M 2^{h j}}\rpa
\lpa \prod_{j= 0}^{m-1} {\pb\lcr |y| \leq d_j\rcr}^{M 2^{h j}}\rpa.$$
To estimate the infinite product, we use the following well-known tail
behavior of $y$ (see e.g. Property 1.2.15 in \cite{ST}):
\begin{equation}
\label{tail}
\lim_{r\uparrow +\infty} r^{\alpha} \pb\lcr |y| > r\rcr\; =\; K_1\,
 \in\, (0, +\infty).
\end{equation}
This yields
\begin{eqnarray*}
\log \prod_{j\geq m} {\pb\lcr |y| \leq d_j\rcr}^{M 2^{h j}} & =
 & M\;\sum_{j\geq m} 2^{h j} \log\lpa 1 - \pb\lcr |y| >
d_j\rcr\rpa\\
& \geq & -C \sum_{j\geq m} 2^{h j}\, \pb\lcr |y| >
2^{(z - \delta)(j-m)}\rcr\\
& \geq & -C \sum_{j\geq m} 2^{h j} 2^{-\alpha(z
- \delta)(j-m)}\\
& \geq & - C\, 2^{h m},
\end{eqnarray*}
where in the last inequality we used $h -\alpha(z -\delta) < 0$.
 The estimate of the finite product is even simpler. Since $y$ has a
 positive density in the neighbourhood of the origin, we have
$$\lim_{\ee\downarrow 0} \ee^{-1} \pb\lcr |y| \leq \ee\rcr\; =\;
K_2\, \in\, (0, +\infty).$$
This entails
\begin{eqnarray*}
\log \prod_{j = 0}^{m-1} {\pb\lcr |y| \leq d_j\rcr}^{M 2^{h j}} &
= & M\;\sum_{j = 0}^{m-1} 2^{h j} \log\pb\lcr |y| \leq d_j\rcr\\
& \geq & -C \sum_{j = 0}^{m-1} 2^{h j}\lpa 1 \, +\, \log\lpa
2^{(z + \delta)(m-j)}\rpa\rpa\\
& \geq & -C 2^{h m}\sum_{j = 0}^{m-1} (m-j)\, 2^{h (j-m)}\\
& \geq & - C\, 2^{h m},
\end{eqnarray*}
where in the last inequality we used $h > 0$. Putting everything
together now yields
$$\pb\lcr |y_{jn}| \leq d_j, \; 1\leq n \leq M 2^{h j}, \; j\geq 0\rcr \; \geq \; \exp - C 2^{h m}$$
for a constant $C$ not depending on $m$.

\fin

\begin{remark} {\em In the Non-Gaussian case, it is easy to see that the
condition $\delta< z-  h/\alpha$ is also necessary,
because of the heavy tails of $\alpha$-stable random variables: indeed,
if the $y_{jn}$'s are mutually independent and if $z = \delta +
h/\alpha$, then it follows from (\ref{tail}) that
$$\pb\lcr |y_{jn}| \leq d_j, \; 1\leq n \leq M 2^{h j}, \; j\geq 0\rcr \; = \; 0.$$
}
\end{remark}

\subsection{Some elements of wavelet theory}
\label{wavelets}
The proof of Theorem \ref{t41} relies mainly on a suitable
wavelet decomposition of $R$, which we recall here for the sake of
completeness. In this paragraph we fix once and for all a semi-norm
$\|.\|\in\LL(\beta, p)$ which is $\ell$-smooth-finite, and $R$ a
continuous RLP satisfying the assumptions of Theorem \ref{t41}.

There exist two functions
$\varphi, \psi\in\CC^{\ell}_K$ ("wavelet parents") such that $\psi$
has vanishing moments up to order $\ell$:
\begin{equation}
\label{mom}
\int_{-\infty}^{\infty} t^k \psi(t)\, dt\; =0\; \qquad 0\le k \le  \ell,
\end{equation}
and such that the wavelet functions
$\lacc \psi_{jn},\; n\in\z,\; j>0\racc$ and
 $\lacc \varphi_{n}, \; n\in\z\racc$, respectively defined by
$$\psi_{jn}(t)\;= \; 2^{j/2}\psi(2^jt-n)\;\;\;\;\;
\mbox{and}\;\;\;\;\;\varphi_{n}(t)\; =\;\varphi(t-n),$$
form an orthogonal base of $L_2(\R)$. We refer to Daubechies' construction
\cite{Dau}, Section 6.4 for the definition and Section 7.1 for smoothness
properties of these compactly supported wavelets. A useful book is also
\cite{M}.

Fix $[-D,D]$ an interval containing the support of $\psi$. Consider $I,$
the integration operator on compactly supported functions:
$$If(t)\;= \;\int_{-\infty}^t f(s)\, ds$$
for every such $f$ and $t \in \R$. We set $I^0$ for the identity
operator and $I^k$ for the $k$-th iteration of $I$, $k\geq 1$.
Since for every $k\geq 0$
$$I^{k+1}f(t)\;= \;\frac{1}{k!}\int_{-\infty}^t f(s) (t-s)^k\, ds,$$
the moment condition (\ref{mom}) on $\psi$
entails that the functions $I^k\psi$ are also supported by the interval
 $[-D,D]$ for $0\le k\le \ell$. In particular,
\begin{equation}
\label{border}
I^k\psi(\pm D)\; = \;0, \qquad
0\le k\le \ell.
\end{equation}
It is well-known that under our assumptions on $H$, the process
$\Un_{[0,D+1]} R$
belongs a.s. to $L_2(\R)$ (see \cite{ST} Chapter 11). Hence, we can write
its wavelet decomposition:
\begin{equation}
\label{ltwo}
\Un_{[0,D+1]} R\; =\; \sum_{j,n} r_{jn} \psi_{jn} +   \sum_{n} r_n
\varphi_{n},
\end{equation}
where
$$ r_{jn}\;=\; \int_0^{D+1} R(s) \psi_{jn}(s)\, ds \;\;\;\;\;\mbox{and}
\;\;\;\;\; r_n\;=\; \int_0^{D+1} R(s) \varphi_{n}(s) ds. $$
Actually more can be said about the convergence of the series on
the right-hand side in (\ref{ltwo}). Namely, again from our assumption on
$H$, it is well-known that $R$ is locally $\eta$-H\"older for some
$\eta > 0$ \cite{KM} \cite{Tak}, and in particular the function
$R$ coincides on $[-D,2D+1]$ with an $L_2$-function which is
{\it globally} $\eta$-H\"older.
We know then (see e.g. Theorem 7 in \cite{M}, Chapter 6) that the series
(\ref{ltwo}) converges to $R$ {\em uniformly} on $[0,1]$.

Besides, in (\ref{ltwo}) we can delete each $\psi_{jn}$ and $\varphi_n$
whose support does not overlap with $[0,1]$, and the remaining series
still converges to $R$ uniformly on $[0,1]$. Since $\varphi$ and $\psi$
have compact support, this leads to a decomposition of $R$ into a
functional array of exponential size: there exists a constant $M$
depending only on $D$ such that
$$R\; =\; \sum_{j\geq 1}\lpa\sum_{|n|\leq M2^j} r_{jn}
\psi_{jn}\rpa\; + \;\sum_{|n|\leq M} r_n \varphi_{n}.
$$

\subsection{Proof of Theorem \ref{t41}}
We may (and will) suppose that $\ell > H$.
For every $n\in\z$ we will set $\psi_{0n} = \phi_n$ and   $r_{0n} = r_n$
for the sake of
concision. By lower semi-continuity and triangle's inequality,
we clearly have
\begin{eqnarray*}
||R|| & \leq & \sum_{j\geq 0}\lva\lva \sum_{|n|\leq M2^j} r_{jn}
\psi_{jn}\rva\rva.
\end{eqnarray*}
Besides, since the parent functions have support in $[-D,D]$,
for each $j\geq 0$
we can split the family
$$\lacc \psi_{jn}, \; |n|\leq M2^j \racc$$
into $2D$ subfamilies such that in each subfamily the functions have
 supports with disjoint interiors. Suppose first that $p < +\infty$.
 Condition (G) yields
\begin{equation}
\label{ineqr}
||R|| \; \leq \; C \,\sum_{j\geq 0}{\lpa\sum_{|n|\leq M2^j}
 |r_{jn}|^p ||\psi_{jn}||^p\rpa}^{1/p}.
\end{equation}
On the one hand, by homogeneity, $\beta$-self-similarity, contractivity
and smooth-finiteness of $||\cdot||$, it is easy to see that for every
$j,n$,
\begin{equation}
\label{smoopsi}
||\psi_{jn}||\; \leq \; C\, 2^{(1/2 + \beta)j}
\end{equation}
for some constant $C$ independent of $j, n$. On the other hand,
 if $j\geq 1$,
\begin{eqnarray}\label{rjn}
r_{jn} &=& 2^{j/2}\int_0^{+\infty}\!\!\!\! R(s)\psi(2^j s - n)\, ds
\nonumber\\
&=& 2^{-j/2} \int_0^{+\infty}\!\!\!\! R(2^{-j} s)\psi(s - n)\,
 ds\nonumber\\
&\elaw& 2^{-(H + 1/2)j} \int_0^{+\infty}\!\!\!\! R(s)\psi(s - n)\,
ds\nonumber\\
& = & 2^{-(H + 1/2)j} \int_{\R} R(s + n)\psi(s)\, ds,
\end{eqnarray}
where in the penultimate inequality we appealed to the $H$-self-similarity
 of $R$. Plugging (\ref{smoopsi}) and (\ref{rjn}) into
(\ref{ineqr}) yields
\begin{equation}
\label{ineqr2}
||R|| \; \leq \; C \,\sum_{j\geq 0} 2^{-(H-\beta-1/p)j}\lpa
\sup_{|n|\leq M2^j} |r_{jn}'|\rpa,
\end{equation}
where $r_{jn}'$ stand for the {\em renormalized} wavelet coefficients of
$R$, viz.
\begin{equation}
\label{reno}
r_{0n}'\; =\; r_{0n} \;\;\;\;\mbox{and}\;\;\;\;
 r_{jn}'\;\elaw\;\int_{\R} R(s + n)\psi(s)\, ds\;\;\mbox{if $j\geq 1.$}
\end{equation}
Notice that one can write the wavelet coefficients
from (\ref{reno}) in the following integral form: if $j\ge 1$, then
$r_{jn}' = \tau_{n} y_{jn}$,
where
\begin{eqnarray*}
\tau_{n}^{\alpha} & = & \int_0^{+\infty} {\lpa \int_{\R} {\lpa u + n - s
\rpa}_+^{H'}\psi(u) \; du\rpa}^{\alpha}\; ds
\end{eqnarray*}
$H'=H-1/\alpha$,
and $\lacc y_{jn}, \; |n| \leq M 2^j, \; j\geq 0\racc$ is an array of
identically distributed symmetric $\alpha$-stable random variables.
 We first aim to prove
\begin{equation} \label{sup}
\sup_{n} \tau_{n}\; < \; +\infty.
\end{equation}
To get this uniform bound, we first recall that $\psi$ has its
support in $[-D, D]$, so that the integral defining
$\tau_{n}^{\alpha}$ can be rewritten as follows:
\begin{eqnarray*}
\tau_{n}^{\alpha} & = & \int_0^{n + D} {\lpa \int_{-D}^D {\lpa u + n - s
\rpa}_+^{H'}\psi(u) \; du\rpa}^{\alpha}\; ds.
\end{eqnarray*}
We cut the domain of integration over $s$ into $[0, n-2D]$ and
$[n-2D, n +D]$. The first integral is given by
$$I^1_{n}\; =\; \int_0^{n - 2D} {\lpa \int_{-D}^D {\lpa u + n - s
\rpa}^{H'}\psi(u) \; du\rpa}^{\alpha}\; ds.$$
We first transform
$$\int_{-D}^D {\lpa u + n - s \rpa}^{H'}\psi (u)\, du$$
through $\ell$ successive integrations by parts. Recalling (\ref{border}),
 we see that each time the border terms
$I^m\psi(\pm D)$, $0\leq m\leq \ell$, vanish. In the end, we obtain
\begin{eqnarray}\label{ijn1}
I^1_{n}& = & C\,\int_0^{n - 2D} {\lpa \int_{-D}^D {\lpa u + n - s \rpa}^{H'
 - \ell}I^{\ell}\psi(u) \; du\rpa}^{\alpha}\; ds\nonumber\\
& \leq & C\,\int_0^{n - 2D} {\lpa D + n - s \rpa}^{\alpha(H' - \ell)}\;
ds\nonumber\\
& \leq & C\,\int_D^{+\infty} {s}^{\alpha(H - \ell) -1}\; ds,
\end{eqnarray}
where in the first inequality we used $\ell>H> H'$.
The integral on the right-hand side is clearly finite and independent of
$n$.

The integral over the second domain is given, after a change of variable, by
\begin{eqnarray}\label{ijn2}
I^2_{n}= \int_0^{3D} {\lpa \int_{-D}^D {\lpa u + 2D - s\rpa}_+^{H'}
\psi(u) \;  du\rpa}^{\alpha} ds &
\leq & C\,  \int_0^{3D} {\lpa \int_0^{3D}u^{H'}\, du\rpa}^{\alpha} ds.
\end{eqnarray}
Since $H' > -1/2$, the integral on the right-hand side is again finite
and independent of $n$. Putting (\ref{ijn1}) and (\ref{ijn2})
together yields (\ref{sup}) as desired.

The following upper bound on
$||R||$ is a direct consequence of (\ref{ineqr2}) and (\ref{sup}):
\begin{equation}
\label{iner}
||R||\; \leq\; C\sum_{j\geq 0} 2^{-(H-\beta - 1/p)j}\lpa\sup_{1\leq
 |n|\leq M 2^j} |y_{jn}|\rpa,
\end{equation}
where $\lacc y_{jn}, \; |n| \leq M 2^j, \; j\geq 0\racc$ is an array of
identically distributed symmetric $\alpha$-stable random variables on $\R$.

The end of the proof is now standard and follows \cite{St2}, Theorem 3.1.
Choose $\delta > 0$ such that
$$\lacc\begin{array}{ll}
\delta <H-\beta-1/p  & \mbox{if $\alpha=2$,}\\
\delta <H-\beta-1/p - 1/\alpha & \mbox{if $\alpha < 2$.}
\end{array} \right.
$$
 Let $m$ be a positive  integer. Set
\begin{equation}
\label{dj}
d_j\; =\; 2^{(H-\beta-1/p) (j-m) - \delta |j - m|}
\end{equation}
for every $j\geq 0$. On the one hand, it is clear from (\ref{iner})
and (\ref{dj}) that
$$\lacc |y_{jn}| \leq d_j, \; 1\leq |n|\leq M 2^j,\; j\geq 0\racc\;
\subset\; \lacc ||R||\, \le \,C(\delta)\, 2^{-(H-\beta-1/p)m} \racc.$$
On the other hand, it follows from Lemma \ref{samo} (with $h = 1$) that
$$\pb\lcr |y_{jn}| \leq d_j, \; 1\leq |n|\leq M 2^j,\; j\geq 0\rcr\;
\geq\;\exp - C 2^m .$$
Since $H-\beta-1/p=1/\gamma$, we obtain
$$\liminf_{m\to \infty} 2^{-m}\log\pb\lcr ||R||\leq C(\delta)
2^{-m/\gamma} \rcr \; > \; -\infty$$
which is equivalent to
$$\liminf_{\ee\downarrow 0}\ee^{\gamma}\log\pb\lcr ||R||\leq \ee \rcr
 \; > \; -\infty,$$
as desired. The proof is complete in the case $p < \infty$.
The case $p = \infty$ can be handled exactly in the same way, replacing
(\ref{ineqr}) by 
\begin{eqnarray*}
||R|| & \leq & C \,\sum_{j\geq 0} \lpa \sup_{|n|\leq M2^j}|r_{jn}| ||\psi_{jn}||\rpa.
\end{eqnarray*}
\fin

\begin{remarks}{\em (a) It is clear from the above proof that
the assumption $\|.\|\in\LL(\beta, p)$ is not necessary to get
the lower bound. We just need $\|.\|\in\Lde(\beta, p)$ w.r.t. $\Psi$,
a wavelet family generated by parents $\{\varphi, \psi\}$ that are smooth
enough. But the required smoothness depends
on the parameter $H$, and this would lead to much heavier notations.
For this reason we prefer stating Theorem \ref{t41} in this form,
save for the loss of generality.
\vspace{2mm}

\noindent
(b) There is at least an example where $\Lde(\beta, p)$
is really more relevant than $\LL(\beta, p)$. It is well-known (see \cite{M}
 Chapter 6) that $(\aa, p, q)$-Besov semi-norms, $\aa > 0$ and $p, q
\geq 1$, are equivalent to sequential norms on the wavelet
coefficients. More precisely,
$$ \left\| \sum_n x_n \phi_{n}+
 \sum_{j,n} x_{jn} \psi_{jn}\right\| \; \sim\;
 \left( \sum_n |x_n|^p \right)^{1/p}
\; + \; \left\| \lacc 2^{(\aa-1/p)j} \left(
 \sum_{n} |x_{jn}|^p
\right)^{1/p}\racc_{j\ge 1}
 \right\|_{\ell_q}.$$
Hence, with our notations,  it is clear that the $(\aa, p, q)$-Besov
semi-norm belongs to $\Lde(\aa -1/p, p)$ w.r.t. any wavelet basis $\Psi$.
Hence, when $p = +\infty$, Theorem \ref{t31} yields the existence of the constant 
for $R$ with $\gamma \;=\; (H-\aa)^{-1}$. Note that this rate is in accordance
with the results of \cite{St2}, which covered the range of
parameters $\alpha=2$ and $0<H<1$ (and with no condition on $p$).

\vspace{2mm}

\noindent
(c) In the Gaussian case, the smooth-finiteness of the semi-norm is obtained for free
when $\|R\|$ is a.s. finite. Indeed, if $\H$ denotes the reproducing kernel
Hilbert space associated with $R$, then it is well-known that $\H$ contains
 functions with $(H+1/2)$-th derivative in $L_2$. On the other hand, from a
0-1 law for Gaussian measures on linear spaces, the finiteness of $\|R\|$
with positive probability yields that $\|f\|<\infty$ for every $f\in \H$.
Hence, $\|.\|$ is $\ell$-smooth-finite as soon as $\ell>H+1/2$. In the
non-Gaussian case, the challenging question whether $\|R\|<\infty$ a.s.
implies the smooth-finiteness of the semi-norm remains open.

\vspace{2mm}

\noindent
(d) For the first efficient use of wavelet methods in similar problems to small
deviation probability, we refer to \cite{AT} where the optimal finite-dimensional approximation of fractional Brownian motion is considered.}

\end{remarks}

\subsection{Proof of Theorem \ref{t42}}

The outline of the proof is the same as that for Theorem
\ref{t41}, except that we use Schauder system
 and provide new estimates for the corresponding coefficients.
Recall that from our assumptions, $0 < H < 2 $ in the Gaussian
case and $1/\alpha<H<2$ in the non-Gaussian case. Again we set
$H' = H - 1/\alpha$. We exclude the case $H' = 0$
where $R$ is Brownian motion: the result follows then directly
from Theorem 1 in \cite{St1}.

Since $R$ is continuous, we can decompose it
on $[0,1]$ along the Schauder system:
$$R_t\; =\; \sum_{j\geq 0} \RC_j(t) + R_1 t,$$
for every $t\in [0,1]$. Here, for every $j\geq 0$,
$$\RC_j(t)\; =\; \sum_{n=1}^{2^j} r_{jn}\,\psi_{jn}(t)$$
where we set $\tk = (n-1)2^{-j}$,
 $\tkk = (n-1/2)2^{-j}$, $\tkkk =n2^{-j}$ and
$$r_{jn} \; =\; 2 R_{\tkk} - R_{\tk} - R_{\tkkk}.$$
We first suppose that $p<\infty$. Condition ($\rm \tilde G$) entails
\begin{equation}
\label{inr}
||\RC_j||\; \leq\; C \ 2^{j\beta}
{\lpa \sum_{n=1}^{2^j}{|r_{jn}|}^p\rpa}^{1/p}.
\end{equation}
Notice that the coefficients $r_{jn}$ can be rewritten as
 $r_{jn}= \sigma_{jn} y_{jn}$, where
\begin{eqnarray} \label{sjnr}
\sigma_{jn}^{\alpha} & = & \int_0^{+\infty} {\lva 2
 {\lpa \tkk - s\rpa}_+^{H'} - {\lpa\tk - s\rpa}_+^{H'}
 - {\lpa\tkkk - s \rpa}_+^{H'} \rva}^{\alpha}\; ds
\end{eqnarray}
(we use the notation $u_+ = u\wedge 0$ for every
 $u\in\R$), and $\lacc y_{jn},\; 1\leq n \leq 2^j, \;
 j\geq 0\racc$ is an array of identically distributed
symmetric $\alpha$-stable random variables. For each $
 1\leq n \leq 2^j, \; j \geq 0$, we will now give an upper bound on
$\sigma_{jn}$ depending only on $j$.

If $n = 1$, then
\begin{eqnarray*}
\sigma_{j1}^{\alpha} & = & \int_0^{2^{-j}} {\lva 2 {\lpa 2^{-(j+1)}
- s\rpa}_+^{H'} - {\lpa 2^{-j} - s \rpa}^{H'} \rva}^{\alpha}\; ds\;
 \leq \; C\, \int_0^{2^{-j}} s^{\alpha H'}\; ds\; = \; C 2^{-\alpha H j}.
\end{eqnarray*}
If $n > 1$, then we set $\tkl = (n-2)2^{-j}\geq 0$ and cut the
domain of integration in (\ref{sjnr}) into $[0,\tkl]$ and
$[\tkl, +\infty)$. Reasoning as above, we see that the second integral
\begin{eqnarray*}
I^2_{jn} & = & \int_{\tkl}^{\tkkk} {\lva 2 {\lpa \tkk - s\rpa}^{H'}_+
- {\lpa\tk - s\rpa}^{H'}_+ - {\lpa\tkkk - s \rpa}^{H'} \rva}^{\alpha}\; ds,
\end{eqnarray*}
is bounded from above by $C\, 2^{-\alpha H j}$. We estimate the first
 integral as follows:
\begin{eqnarray*}
I^1_{jn} & = & \int_0^{\tkl} {\lva 2 {\lpa \tkk - s\rpa}^{H'} -
{\lpa\tk - s\rpa}^{H'} - {\lpa\tkkk - s \rpa}^{H'} \rva}^{\alpha}\; ds
\\
&=& 2^{-j}   \int_2^{n} {\lva 2 {\lpa \tkk - \tkkk+2^{-j}u\rpa}^{H'} -
{\lpa\tk - \tkkk+2^{-j}u \rpa}^{H'} -
\lpa 2^{-j}u \rpa^{H'} \rva}^{\alpha}\; du
\\
&=& 2^{-\alpha Hj}   \int_2^{n} {\lva 2 {\lpa u-1/2 \rpa}^{H'} -
{\lpa u-1 \rpa}^{H'} -
{u}^{H'} \rva}^{\alpha}\; du
 \\
&\le& C\ 2^{-\alpha Hj}   \int_2^{\infty} u^{(H'-2)\alpha}\; du
 \\
&=& C\ 2^{-\alpha Hj},
\end{eqnarray*}
where the last equality comes from $H<2$.
Putting everything together yields now
$$\sigma_{jn}\; \leq \; C\, 2^{-H j}$$
and, recalling (\ref{inr}),
\begin{equation}
\label{inr2}
||\RC_j||\; \leq\; C 2^{(\beta + 1/p - H)j}\lpa
\sup_{1\leq n\leq 2^j} |y_{jn}|\rpa.
\end{equation}
By lower semi-continuity and triangle's inequality, this entails
\begin{eqnarray*}
||R||\; \leq\; C \lpa |R_1|\; +\; \sum_{j\geq 0}
2^{(\beta + 1/p - H)j}\lpa\sup_{1\leq |n|\leq M 2^j} |y_{jn}|\rpa\rpa
\end{eqnarray*}
and since $R_1$ has a symmetric $\alpha$-stable law,
we can finish the proof as in Theorem \ref{t41}.

\fin

\begin{remarks}{\em (a) In the latter proof, the estimate on
$I^1_{jn}$ becomes too crude when $H\geq 2$. This is the principal
reason why we need to introduce smoother wavelets in Theorem \ref{t41}.

\vspace{2mm}

\noindent
(b) As we said before, the proofs of Theorems \ref{t41} and \ref{t42}
work only for continuous processes. The main reason for
this comes from Lemma \ref{samo}. On the one hand, the exclusion of
the boundary value $z = \delta + h/\alpha$ in the non-Gaussian case
cancels important case $H = 1/\alpha$.
On the other hand, \v Sid\'ak's inequality
makes us handle the $\|\RC_j\|$'s as if they were {\em independent}.
But if $R$ is discontinuous, then it is possible that its jumps have
a significant influence on every level  $\RC_j$, so that our estimate
of the series $\sum_j \|\RC_j\|$ may not be realistic anymore.
\vspace{2mm}

\noindent
(c) In the Non-Gaussian case and when $H < \beta + 1/p + 1/\alpha$,
it is easy to see that the statements of Theorems \ref{t41} and
\ref{t42} are false, even if $H > \beta + 1/p$. Suppose for example
that $H < 1/\alpha$. It is well-known (see \cite{ST} Chapter 10) that a.s.
$$\sup_{0\leq t\leq 1} |R_t|\; = \; +\infty,$$
so that Theorems \ref{t41} and \ref{t42} cannot hold with respect
to the supremum norm (i.e.
$\k = +\infty$ in Theorem \ref{t31}), although here $\beta = 0$ and
$p = +\infty$. Similarly, if $H = 1/\alpha$, then clearly
$$\sup_{0\leq s <t\leq 1}\lpa \frac{|R_t - R_s|}{{|t-s|}^{\eta}}\rpa\;
= \; +\infty\;\;\;\;\mbox{a.s.}$$
for every $\eta > 0$, so that Theorems \ref{t41} and \ref{t42} cannot
hold with respect to any $\eta$-H\"older semi-norm ($\beta =\eta,
p = +\infty$). It is natural to conjecture that these two theorems
remain true in full generality when $H = \beta + 1/p + 1/\alpha$,
but in view of the above Remark (b), we probably need different methods.}
\end{remarks}

\section{Lower bounds for LMP: Finiteness of the constant for continuous
LFSM}  \label{motion}

Our aim in this section is to extend Theorem \ref{t42}
to $X$, the continuous LFSM which we defined in Section 2.
Recall that $X$ admits an independent
decomposition $X = R + M$ where $R$ (resp. $M$) is an RLP (resp. an LMP)
 with the same parameters. The following theorem can be viewed as
 a generalization of Lemma 3.2 in \cite{LL1}.

\begin{theorem} \label{t51} Let $||\cdot||\in\Lde(\beta, p)$ w.r.t.
the Schauder system and $M$ be an LMP with parameters $\alpha\in (1, 2]$,
$H\in (1/\alpha,1)$. Suppose that $H > \beta + 1/p$ if $\alpha = 2$ and
$H > \beta+ 1/p+1/\alpha$ if $\alpha\neq 2$. Then
$$\lim_{\ee\downarrow 0}\ee^{\gamma}\log\pb\lcr ||M||\leq \ee \rcr
 \; = \; 0,$$
where $\gamma=(H-\beta-1/p)^{-1}$.
\end{theorem}

Using Theorems \ref{t31} and \ref{t42} along with the elementary
independence argument developed in \cite{LL1} p. 1334,
Theorem \ref{t51} yields readily the desired small deviation
theorem for $X$:

\begin{theorem} \label{t52} Let $||\cdot||\in\Nde(\beta, p)$ w.r.t.
the Schauder system and $X$ be a continuous LFSM with parameters
$\alpha\in (1, 2]$, $H\in (1/\alpha,1)$. Suppose that $H > \beta + 1/p$
if $\alpha = 2$ and $H > \beta+ 1/p+1/\alpha$ if $\alpha\neq 2$. Then
$$\lim_{\ee\downarrow 0}\ee^{\gamma}\log\pb\lcr ||X||\leq
\ee \rcr
 \; = \; \lim_{\ee\downarrow 0}\ee^{\gamma}\log\pb\lcr ||R||\leq \ee \rcr
 \; = \; -\k$$
where $\k\in (-\infty, 0),
\gamma = (H-\beta-1/p)^{-1}$ and $R$ is the RLP associated to $X$.
\end{theorem}

We now proceed to the proof of Theorem \ref{t51}, which is quite
analogous to the proof of Theorem \ref{t42}, except that our
estimates on the Schauder coefficients will not be uniform in
time-argument since $M$ has a singularity at 0.

\subsection{Proof of Theorem \ref{t51}}

Again we set $H' = H- 1/\alpha\in (-1,1)$ and exclude the trivial
situation when $H' = 0$. Since $M$ is continuous, we can decompose
it on $[0,1]$ along the Schauder system: for every $t\in [0,1]$
$$M_t\; =\; \sum_{j\geq 0}  \MC_j(t) + M_1 t,$$
where
$$\MC_j(t) = \sum_{n=1}^{2^j} m_{jn}\,\psi_{jn}(t)$$
and $m_{jn} \; =\; 2 M_{\tkk} - M_{\tk} - M_{\tkkk}$ as in the proof
of Theorem \ref{t42}. We first consider the case $p < +\infty$.
Condition ($\rm \tilde G$) entails
\begin{equation} \label{inm}
||\MC_j||\; \leq\; C 2^{\beta j}{\lpa
\sum_{n=1}^{2^j}{|m_{jn}|}^p\rpa}^{1/p}
\end{equation}
and we just need to evaluate the coefficients
$m_{jn} = \sigma_{jn} y_{jn}$, where
\begin{eqnarray*}
\sigma_{jn}^{\alpha} & = & \int_0^{+\infty} {\lva 2 {\lpa s +
 \tkk\rpa}^{H'} - {\lpa s + \tk\rpa}^{H'} -
 {\lpa s + \tkkk\rpa}^{H'} \rva}^{\alpha}\; ds
\\
&=& 2^{-H\alpha j} \int_0^{+\infty} {\lva 2 {\lpa u + n- 1/2 \rpa}^{H'}
- {\lpa u+n-1\rpa}^{H'} -  {\lpa u + n\rpa}^{H'} \rva}^{\alpha}\; du,
\end{eqnarray*}
and $\lacc y_{jn},\; 1\leq n \leq 2^j, \; j\geq 0\racc$ is an array
of identically distributed symmetric $\alpha$-stable random
variables. Notice first that from the so-called Cooper's formula
$$2 {\lpa u + n- 1/2 \rpa}^{H'} - {\lpa u+n-1\rpa}^{H'} -
 {\lpa u + n\rpa}^{H'} =
- H'(H'-1)
\int_0^{1/2} \int_{u+n-1+\theta}^{u+n-1/2+\theta} v^{H'-2} dv d\theta$$
Clearly, since $H' < 2$, this entails
$${\lva 2 {\lpa u + n- 1/2 \rpa}^{H'} - {\lpa u+n-1\rpa}^{H'} -
 {\lpa u + n\rpa}^{H'} \rva} \le
C \min \left\{ (u+n-1)^{H'-2} ,1 \right\}.$$
We now fix some $h \in (0,1)$ and first estimate $\sigma_{jn}$ when $2^{h j} < n \leq 2^j$:
\begin{eqnarray}
\sigma_{jn}^{\alpha} & \le &  C 2^{-H\alpha j}
\int_0^{+\infty} (u+n-1)^{\alpha(H'-2)}  \; du
\nonumber \\
&\le&   C 2^{-H\alpha j} (n-1)^{\alpha(H'-2)+1}
\nonumber \\
&=&   C 2^{-H\alpha j} 2^{-h\alpha(2-H)j}\ .
\label{sig1}\end{eqnarray}
The estimate of $\sigma_{jn}$ when $1 \leq n\leq 2^{h j}$
is even simpler:
\begin{eqnarray}
\sigma_{jn}^{\alpha} & \le &  C 2^{-H\alpha j}
\lpa
\int_1^{+\infty} u^{\alpha(H'-2)}  \; du
+
\int_0^{1} C  \; du \rpa\;\le\; C 2^{-H\alpha j} .
\label{sig2}
\end{eqnarray}
Recalling (\ref{inm}), we get
\begin{eqnarray*}
||\MC_j||  &\leq&  C \, 2^{\beta j} {\lpa
\sum_{n=1}^{2^{h j}} \sigma_{jn}^p   +
\sum_{n = 2^{h j}+1}^{2^j} \sigma_{jn}^p
\rpa}^{1/p}  \sup_{1\le n\le 2^j} |y_{jn}|
\\
&\le& C\, 2^{-(H-\beta) j} {\lpa
2^{h j}   + 2^{j} 2^{-h(2-H)pj}
\rpa}^{1/p}  \sup_{1\le n \le 2^j} |y_{jn}|.
\end{eqnarray*}
The balance is attained at $h=(1+(2-H)p)^{-1}\in (0,1)$, whence
\begin{equation}\label{inm2}
||\MC_j||  \leq
 C \, 2^{-(H-\beta-h/p) j}  \sup_{1\le n \le 2^j} |y_{jn}|.
\end{equation}
Using (\ref{inm2}) instead of (\ref{inr2}), we can proceed as in
the proof of Theorem \ref{t42} and obtain
$$\liminf_{\ee\downarrow 0}\ee^{\gamma'}\log\pb\lcr ||M||\leq \ee \rcr
 \; > \; -\infty,$$
with $\gamma'=(H-\beta-h/p)^{-1}< \gamma$. In particular
$$\liminf_{\ee\downarrow 0}\ee^{\gamma}\log\pb\lcr ||M||\leq \ee \rcr
 \; = \; 0,$$
which completes the proof of Theorem \ref{t51} for $p<\infty$.

\vspace{2mm}

The case $p = +\infty$ requires more careful estimates.
Condition ($\rm \tilde G$) entails
\begin{equation} \label{inm5}
||\MC_j||\; \leq\; C 2^{\beta j}\,\sup_{1\leq n\leq 2^j} |m_{jn}|.
\end{equation}
Again we fix some $h\in(0,1)$. Using (\ref{inm5}), (\ref{sig1}),
 and (\ref{sig2}) we can write,
\begin{eqnarray*}
 ||\MC_j||\; &\leq&\;  C\, 2^{\beta j}
\sup_{1\le n \le 2^j} \sigma_{jn} |y_{jn}|
\\
&\leq&   C\, 2^{-(H-\beta) j}  \lpa
\sup_{1\le n \le 2^{hj}}  |y_{jn}|  +
2^{-h_1 j}\sup_{2^{hj} < n \le 2^j} |y_{jn}|  \rpa,
\end{eqnarray*}
where we set $h_1=(2-H)h$. For every integer $m>0$ we focus on the event
$$\Omega_m\; = \;\lacc |y_{jn}| \leq d_{jn}(m), \; 1\leq n \leq
2^j, \; j \geq 0\racc$$
where we set $m_1=[m/h]$, $\delta\in (0, H-\beta)$ and
$$d_{jn}(m)\; =\; \lacc\begin{array}{ll}
2^{(H-\beta)(j-m_1)-\delta|j-m_1|}& \mbox{if $1\leq n \leq 2^{hj}$}\\
2^{(H-\beta+h_1)(j-m)-\delta|j-m|}& \mbox{if $2^{hj} < n \leq 2^{j}.$}
                              \end{array}\right.$$
Take now $h$ small enough such that
$$H-\beta < H- \beta + h_1 < (H -\beta)/h.$$
On the one hand, it is clear that
\begin{equation} \label{inm3}
\Omega_m\; \subset\; \lacc \sum_{j\ge 0} ||\MC_j||\; \leq\;
C 2^{-(H-\beta+h_1)m} \racc.
\end{equation}
On the other hand, it follows from decorrelation argument of
Lemma \ref{samo} that
\begin{equation} \label{inm4}
\log\pb\lcr \Omega_m\rcr\; \geq
- \, C\, \lpa  2^{h m_1} + 2^{m} \rpa\; \geq\; -\, C \, 2^m.
\end{equation}
Using (\ref{inm3}) and (\ref{inm4}), we can now finish the proof
exactly as in the case $p < +\infty$.

\fin

\begin{remarks}{\em (a) Contrary to \cite{LL1}, the above proof
does not require any entropy argument and relies only on an elementary
estimate of the Schauder coefficients.

\vspace{2mm}

(b) It would be quite interesting to calculate the optimal rate
$$\gamma_0\; =\; \inf\lacc\gamma > 0\;\lva\; \lim_{\ee\downarrow
0}\ee^{\gamma}\log\pb\lcr ||M|| \leq \ee \rcr\; = \; 0\racc.\right.$$
We did not study this question in details but it seems plausible that
$\gamma_0 = 0$: indeed, $M$ is a $\CC^{\infty}$ process on $(0, +\infty)$,
and one may expect subexponential rate for its small deviation
probabilities.}
\end{remarks}

\section{Examples: particular semi-norms} \label{survey}

In this section, we place all our results in the context of the previous
literature, and show what is new. We try to be as exhaustive as possible
 as far as RLP's and LFSM's are concerned. We refer to the surveys \cite{LS}
 \cite{Lf} for more information about other processes and for further
references.

Everywhere $R$ (resp. $X$) will be an RLP (resp. an LFSM) with Hurst
parameter $H > 0$ and stability index $\alpha\in (0,2]$, whereas
$\|.\|$ will be a semi-norm in $\NN(\beta, p)$ or $\Nde(\beta, p)$
for some $\beta\in\R$ and $p\in (0, +\infty]$. We always suppose that
$H > \beta + 1/p$ and set
$$\gamma = (H-\beta-1/p)^{-1}$$
for our small deviation rate. $\k$ will stand for the small deviation
constant appearing on the right-hand side of Theorem \ref{t31}.
Sometimes we will call them just "the rate" and "the constant".
Let us begin with the most classical case, which deserves of course
a particular mentioning.

\subsection{Brownian motion}

Brownian motion $B$ is an RLP (or an LFSM) with parameters $H=1/2$
and $\alpha=2$, hence it clearly satisfies the assumptions of
Corollary \ref{c42} (or Theorem \ref{t52}). Notice that here our
Theorem \ref{t42} just amounts to Stolz's lower bound criterion
(Theorem 1 in \cite{St1}), so that in the present paper the novelty
comes only from Theorem \ref{t31}. Corollary \ref{c42} entails
the existence of Brownian small deviation constants under almost
all the
classical semi-norms (with the regrettable exception of certain Besov semi-norms
for which we were unable to prove the appropriate superadditivity index).
In other words, we get
$$\lim_{\ee\downarrow 0}\ee^{\gamma}\log\pb\lcr ||B||\leq \ee \rcr
 \; =\; - \k\;\in\; (-\infty, 0),$$
where the dependence of $\gamma$ on the semi-norm $||\cdot||$ is given
by the following table:

\vspace{2mm}

\begin{center}
\begin{tabular}{|c|c|}
\hline
$||\cdot||$ & $\gamma$ \\
\hline
\hline
Supremum & 2 \\
\hline
$L_p$ & 2 \\
\hline
$\aa$-H\"older & $2/(1- 2\aa)$ \\
\hline
$p$-variation & $2p/(p- 2)$ \\
\hline
$(\aa, p)$-Sobolev & $2/(1- 2\aa)$ \\
\hline
$(\aa, \infty, q)$-Besov & $2/(1- 2\aa)$ \\
\hline
\end{tabular}
\end{center}

\vspace{2mm}

The existence of a finite small deviation constant for
Brownian motion under $L_p$ and supremum norms is a very classical
result, and in this situation the constant $\k$ can even be more or less
explicitly computed (see \cite{LS} and the references therein).
The existence of the constant under H\"older semi-norms was first
proved by Baldi and Roynette \cite{BR}, and it is actually a direct
consequence
of Stolz's upper bound criterion (Theorem 2 in \cite{St1}) and of the
standard subadditivity argument involving the Markov property of
 Brownian motion. In the three last examples, Stolz's results do provide
the right rate, viz.
$$ -\infty \; < \; \liminf_{\ee\downarrow 0}\ee^{\gamma}\log\pb\lcr
||B||\leq \ee \rcr \; \leq \; \limsup_{\ee\downarrow 0}
\ee^{\gamma}\log\pb\lcr ||B||\leq \ee \rcr \; < \; 0$$
with the above corresponding $\gamma$, but here the existence
of the limit does not follow directly. In these three cases,
our result seems to be new.

It is of course quite tantalizing to {\em compute} these new
Brownian constants. As a rule, this kind of (hard) problems
requires analytical techniques depending heavily on the choice
of the semi-norm, and it seems difficult to tackle them through
a global study.

We stress that the above table matches quite accurately the sample
 path properties of Brownian motion. It is well-known that $B$
is $\aa$-H\"older if and only if $\aa < 1/2$ and has finite $p$-variation
 if and only if $p > 2$. Besides, it has finite
$(\aa, p)$-Sobolev semi-norm if and only if $\aa < 1/2$
\cite{CHS}.

Besov semi-norms deserve a special remark since no necessary
and sufficient conditions seem to be available. Roynette \cite{Ro}
had proved that Brownian motion has finite $(\aa, p, q)$-Besov
semi-norm whenever one of the three following situations occurs:

\vspace{2mm}

\begin{center}
\begin{tabular}{|c|c|c|}
\hline
$\aa$ & $p$ & $q$ \\
\hline
\hline
$<1/2$ & $>1/\aa$ & $\ge 1$\\
\hline
$<1/2$ & $=1/\aa$ & $=1$ \\
\hline
$=1/2$ & $>2$ & $ = +\infty$ \\
\hline
\end{tabular}
\end{center}
\vspace{2mm}

\noindent
Notice that the latter situation is rather specific,
 since here small deviation probabilities vanish when $\varepsilon$ is
small enough, thus the investigation is pointless.

We are able to prove the existence of a finite small deviation constant in the following situation:

\vspace{2mm}

\begin{center}
\begin{tabular}{|c|c|c|}
\hline
$\aa$ & $p$ & $q$ \\
\hline
\hline
$<1/2$ & $= +\infty$ & $\geq 1$ \\
\hline
\end{tabular}
\end{center}

\subsection{Gaussian fractional processes}

In this paragraph $\alpha = 2$ and $H > 0$, the case $H =1/2$ being
implicitly excluded. Up to normalization constant our process $X$ is
the usual fractional Brownian motion, while $R$ can be viewed as the
fractionally integrated Brownian motion. Before stating
our results, we will give a list of local properties of $R$,
most of which are classical. As we mentioned in Section 2, all these
local properties remain true for $X$ as soon as $H < 1$.

\vspace{2mm}

$R$ is $\aa$-H\"older (resp. $\aa$-Lipschitz) if and only if $H > \aa$.
In particular, $R$ has finite $p$-variation
if $H> 1/p$ and only if $H\geq 1/p$ (the boundary
value being excluded when $H\geq 1/2$ - see Theorem 5.4 in \cite{DN},
 but the situation $H = 1/p < 1/2$ seems to have escaped investigation).
Finally, $R$ has finite $(\aa,p)$-Sobolev semi-norm if and only if
$H >\aa$ \cite{CHS}.

Less is known about Besov semi-norms. In \cite{CKR}, Proposition 4.1 and
Theorem 4.3, it is proved that if $H < 1$ and
$0 < 1/p < \aa$, then $X$ has finite $(\aa,p,\infty)$-Besov
 semi-norm if and only if $H\geq\eta$ (when $H = \eta$, the small
deviation probabilities vanish for small $\ee$).

\vspace{2mm}

With the same notations as in the preceding paragraph, the following
table is a direct consequence of Corollary \ref{c41} and of Remark (b)
 at the end of paragraph 4.3. Recall that it is available for $R$ and
for $X$ with the same finite constant $\k$, as soon as $H > \beta + 1/p$.

\vspace{2mm}

\begin{center}
\begin{tabular}{|c|c|}
\hline
$||\cdot||$ & $\gamma$ \\
\hline
\hline
Supremum & $1/H$ \\
\hline
$L_p$ & $1/H$ \\
\hline
$\aa$-H\"older & $1/(H-\aa)$ \\
\hline
$p$-variation & $p/(Hp- 1)$ \\
\hline
$(\aa, \infty, q)$-Besov & $1/(H-\aa)$ \\
\hline
\end{tabular}
\end{center}

If in addition $H < 2$, then we can add another line:

\vspace{2mm}

\begin{center}
\begin{tabular}{|c|c|}
\hline
$||\cdot||$ & $\gamma$ \\
\hline
\hline
$(\aa, p)$-Sobolev & $1/(H-\aa)$ \\
\hline
\end{tabular}
\end{center}

\vspace{2mm}

The existence of $\k$ for $X$ under the supremum norm was established
in \cite{LL1} and \cite{S2}. In \cite{LL1}, it is also proved for $R$
when $H < 1$. When $H\geq 1$, the existence of $\k$ is mentioned in
\cite{LL2}. Under $L_p$-norms, the right rate
for $R$ can be found in \cite{BL}. Specialists seem to be aware of
the existence of the constant \cite{LfL2}, but up to our knowledge
this is not explicitly proved in literature. Notice that the case
$p=2$ is rather special because of its Hilbert space background. Quite
recently, the {\em exact} value of
$\k$ under the $L_2$-norm was found for integrated Wiener process
\cite{CL} \cite{KS}, and for fractional Brownian motion \cite{Br}.

\vspace{2mm}

In all other situations, our results seem to be new.

\subsection{Non-Gaussian stable L\'evy processes}

Non-Gaussian symmetric $\alpha$-stable L\'evy processes are RLP's with
parameters $\alpha < 2$ and $H = 1/\alpha$. They have a.s.
discontinuous sample paths, hence in the present paper
they are only concerned by Theorem \ref{t31}.

Recall that they have finite $p$-variation if and only if $p > \alpha$
\cite{BG}. In the case $1<\alpha<2$, it is shown in \cite{CKR} Theorem
6.1 that they have finite $(1/\alpha, p, \infty)$-Besov norm if and only
if $1\leq p <\alpha$.

Our Theorem \ref{t31} yields the following table which can be read as
in the preceding paragraph, except that here we know nothing
{\em a priori} about the finiteness of
$\k$.

\vspace{2mm}

\begin{center}
\begin{tabular}{|c|c|}
\hline
$||\cdot||$ & $\gamma$ \\
\hline
\hline
Supremum & $\alpha$ \\
\hline
$L_p$ & $\alpha$ \\
\hline
$p$-variation & $\alpha p/(p- \alpha)$ \\
\hline
\end{tabular}
\end{center}

\vspace{2mm}

Under the supremum norm, the finiteness of $\k$ dates back to Taylor
\cite{Ta} and Mogulskii \cite{Mo}. Notice that in this case the existence
of
$\k$ follows readily from the Markov property, so that our Theorem
\ref{t31} is a bit useless here. Under $L_p$-norms, the existence
and the finiteness of
$\k$ seems to belong to some ancient folklore,
but we could not find any precise reference in the literature.
We remark that this follows readily from Theorem \ref{t31} and
Taylor-Mogulskii's lower bound. Z. Shi also indicated us that in this
case  the problem of
$L_p$-small deviations can be reduced to some specific
large deviation probabilities for occupation measures which were
previously studied by Donsker and Varadhan. Finally, it is shown in
\cite{Shi} that the case $p=2$ is again rather specific since it can
be derived from the
$L_p$-small deviations of Wiener process,
thanks to an identity in law due to Donati-Martin, Song, and Yor \cite{DSS}.

Nothing is known about the finiteness of $\k$
for the $p$-variation norm, but there is no reason to believe that
the rate $\gamma$ indicated above is not the
right one. This will be the matter of further study.

\subsection{Non-Gaussian fractional processes}

In this paragraph $\alpha < 2$ and $H > 0$, the case $H = 1/\alpha$
being implicitly excluded. As in the Gaussian framework, we begin with
listing local properties of $R$, most of which can be found in
the monograph \cite{ST}. Again, all these local properties remain true
for $X$ when it is well-defined.

\vspace{2mm}

$R$ is sample-bounded iff $H \geq 1/\alpha$ and belongs locally
to $L_p$ if and only if $H > 1/\alpha -1/p$. $R$ is $\aa$-H\"older
(resp. $\aa$-Lipschitz) if and only if $H\geq \aa +1/\alpha$. This
latter fact is less classical and was proved for $X$ by Takashima
\cite{Tak} (see also \cite{KM}). Notice that the ''if" part is quite
straightforward from
$$|R_t - R_s |\; \le\; C\,
||Z||_{\infty}\, |t - s|^{\min\{1,H-1/\alpha\}}$$
for all $s, t\in [0,1]$. It is interesting to note that contrary to
the Gaussian case, the boundary value $\eta=H-1/\alpha$
is included here.
In particular, this entails that $R$ has finite $p$-variation
if and only if $H\geq 1/p +1/\alpha$ \cite{ChG}.

Nothing seems to be known about Sobolev or Besov semi-norms.
From H\"older continuity, we see that $R$ has finite $(\aa, p)$-Sobolev
semi-norm if $H > \aa + 1/\alpha$, yet this does not seem
to be a sharp estimate.

\vspace{2mm}

When $R$ is continuous and when $H > \beta + 1/p + 1/\alpha$, we
get exactly the same tables as in paragraph 6.2. 

\vspace{2mm}

The only
result which we are aware of in this direction is due to
G. Samorodnitsky \cite{Sa}, who had given just a bit less precise
bounds in the small
deviation problem for LFSM's ($1/\alpha<H<1$) under the supremum norm.
 In the present paper we find the right rate,
prove that the constant exists, and our condition on
$H$ is the best possible
since $R$ is not bounded anymore when $H<1/\alpha$. In general,
our results match quite accurately the local sample-path properties
 of $R$. Namely, except for Sobolev and Besov norms, the only
situations which are not covered are the critical values of $H$
w.r.t. H\"older and $p$-variation semi-norms.
\vspace{2mm}

When $R$ is not continuous, i.e. when $H < 1/\alpha$,
the only relevant classical semi-norms are the $L_p$-norms,
as soon as $H>1/\alpha-1/p$. Our Theorem \ref{t31} yields
$\gamma=1/H$ but nothing is known about the finiteness of
the constant $\k$.
We believe that $1/H$ is still the right rate in this
ultimate situation.

\section{Concluding remarks}

(a)  In general, it is not true that $H$ is the only
parameter of a self-similar stable process which is involved in the exponential behavior
of its small deviation probabilities. Indeed, \cite{Sa} delivers an example of an
$\alpha$-stable, $H$-sssi process whose small
deviation rate depends both on $H$ and $\alpha$. Let us also
recall the general results of \cite{Ry} about stable measures,
which entail that
$$\gamma\;\le\; \frac \alpha{1-\alpha}$$
for arbitrary $\alpha$-stable processes with $0<\alpha<1$.

\vspace{2mm}

\noindent
(b) The {\it symmetry} assumption on the stable processes
is a crucial one,
since our proofs are heavily based on two Gaussian results:
Anderson's and \v Sid\'ak's inequalities. Nevertheless, it
seems likely that everything remains true as soon as $Z$ is not
 a subordinator.

\vspace{2mm}

\noindent
(c) Theorem \ref{t31} can be easily extended to
multi-dimensional RLP's having some symmetry,
since they also share the crucial extrapolation-homogeneity
property. However, adaptating the proof of Theorem \ref{t41} in a
multi-dimensional context raises some technical difficulties
related to \v Sid\'ak's inequality.

\vspace{2mm}

\noindent
(d) As we indicated before, there are many different
generalizations of fractional Brownian motion in
non-Gaussian case. If $\tilde X$ is the well-balanced
linear fractional stable motion defined in (\ref{balanced}),
then the proof of Theorem \ref{t31} fails in the absence of
extrapolation-homogeneity. However, we believe that the
rates should be the same.

\vspace{2mm}

\noindent
(e) Although our results are designed for nice translation-invariant
semi-norms, it is worthwhile to mention that there exist elaborated
techniques concerning norms which do not share any translation-invariance.
See especially \cite{LfL1}, \cite{LfL2} and \cite{Shi},
where weighted $L_p$-norms are handled respectively for
Brownian motion, fractional Brownian motion, and symmetric
$\alpha$-stable processes.
\bigskip

\noindent
{\em Acknowledgements:} We are much indebted to Charles Suquet for
valuable advice concerning wavelet techniques, and to Werner Linde
and Zhan Shi
for helpful discussions.  The first-named author was partially
supported by grants RFBR 02-01-00265 and NSh-2258.2003.1.
This work was carried out during a visit at the Technische
Universit\"at Berlin of the second-named author,
who would like to thank Michael Scheutzow for his kind hospitality.

\vspace{10mm}

{\footnotesize
\baselineskip=12pt

\noindent
\begin{tabular}{lll}
& \hskip20pt Mikhail A. Lifshits
    & \hskip40pt Thomas Simon \\
& \hskip20pt Faculty of Mathematics and Mechanics
    & \hskip40pt  \'Equipe d'Analyse et Probabilit\'es\\
& \hskip20pt St-Petersburg State University
    & \hskip40pt  Universit\'e d'\'Evry-Val d'Essonne\\
& \hskip20pt 198504, Stary Peterhof
    & \hskip40pt  Boulevard Fran\c{c}ois Mitterrand\\
& \hskip20pt Bibliotechnaya pl., 2
    & \hskip40pt 91025 \'Evry Cedex\\
& \hskip20pt Russia
    & \hskip40pt France \\
& \hskip20pt E-mail: {\tt lifts@mail.rcom.ru}
    & \hskip40pt E-mail: {\tt simon@maths.univ-evry.fr}
\end{tabular}

}

\end{document}